\documentclass[a4paper,11pt]{article}

\usepackage[OT1]{fontenc}
\usepackage{amsfonts}
\usepackage{amsmath}
\input epsf
\usepackage{epsfig}
\usepackage{latexsym}
\usepackage{graphicx}
\makeatletter
\renewcommand\theequation{\thesection.\arabic{equation}}
\@addtoreset{equation}{section}
\makeatother

\newtheorem{thm}{\bf Theorem}[section]
\newtheorem{lem}{\bf Lemma}[section]

\newtheorem{cor}{\bf Corollary}[section]
\newtheorem{defi}{\bf Definition}[section]
\newtheorem{rem}{\bf Remark}[section]
\newtheorem{hyp}{\bf Hypothesis}[section]

\newcommand{\R}{\mathbb{R}}
\newcommand{\CC}{\mathbb{C}}
\newcommand{\Z}{\mathbb{Z}}
\newcommand{\N}{\mathbb{N}}

\newcommand{\Frac}[2]{{\displaystyle\frac{ #1}{ #2}}} 
\newcommand{\TFrac}[2]{{\textstyle\frac{ #1}{ #2}}} 

\newcommand{\suml}{\sum\limits}

\newcommand{\dindice}[2]{{\stackrel{\scriptstyle #1}{\scriptstyle #2}}}

\newcommand{\supl}{\sup\limits}

\newcommand{\grandO}[1]{\mathop{\hbox{${\cal O}$}}\limits_{#1}}

\newcommand{\I}{{\rm i}}
\newcommand{\E}{{\rm e}}
\newcommand{\eps}{\varepsilon}

\newcommand{\Id}{{\rm Id}}

\newcommand{\Img}{{\rm Im}\;}
\renewcommand{\ker}{\textrm{Ker}\;}

\newcommand{\scal}[2]{\left \langle #1, #2 \right   \rangle}
\newcommand{\Norme}[2]%
  {\left |  #1 \right |_{\hspace{-0.9ex}\raisebox{-0.5ex}{ $\scriptstyle
{#2}$}}}
\newcommand{\NNorme}[2]%
{\left \|  #1 \right \|_{\hspace{-0.9ex}\raisebox{-0.5ex}{ $\scriptstyle
#2  $}}}
\newcommand{\norme}[2]%
  {\left |  #1 \right |_{\hspace{-0.9ex}\raisebox{-0.5ex}{ $\scriptscriptstyle
{#2}$}}}
\newcommand{\nnorme}[2]%
{\left \|  #1 \right \|_{\hspace{-0.9ex}\raisebox{-0.5ex}{ $\scriptscriptstyle
#2  $}}}
\newcommand{\sNorme}[2]%
  { |  #1 |_{\hspace{-0.9ex}\raisebox{-0.5ex}{ $\scriptstyle #2$}}}
\newcommand{\sNNorme}[2]%
{ \|  #1  \|_{\hspace{-0.9ex}\raisebox{-0.5ex}{ $\scriptstyle #2$}}}

\newcommand{\bNorme}[2]%
  { \Bigl |  #1 \Bigr |_{\hspace{-0.9ex}\raisebox{-0.5ex}{ $\scriptstyle #2$}}}
\newcommand{\bNNorme}[2]%
{ \Bigl \|  #1 \Bigr   \|_{\hspace{-0.9ex}\raisebox{-0.5ex}{ $\scriptstyle #2$}}}

\def\V{V}
\def\NN{\mathcal{N}}
\def\Nn{N}
\def\PPhi{\underline{\Phi}}
\def\pphi{\varphi}
\def\m{m}
\def\n{n}
\def\q{q}
\def\l{\ell}
\def\Nu{\Omega}
\def\sig{\sigma}
\def\T{T}
\def\c{c}
\def\rh{\rho}
\def\ga{\gamma} 
\def\gam{\gamma} 
\def\ta{\tau}
\def\be{\beta}
\def\a{a}
\def\k{k}
\def\i{i}
\def\j{j}
\def\p{p}
\def\b{b}
\def\Ezero{E_0}
\def\Eun{E_1}
\def\uzero{u_0}
\def\uun{u_1}
\def\vun{v_1}
\def\Mzero{M_0}
\def\M{M}
\def\y{y}
\def\X{X}
\def\F{F}
\def\om{\omega}

\def\Nrond{\mathcal{N}}
\def\Reste{R}
\def\CSob{\mathcal{C}}
\def\C{C}
\newcommand{\cqfd}{\hfill$\Box$}
\newcommand{\Sob}{\textit{H}}
\newcommand{\A}{\mathcal{A}}
\newcommand{\AAL}{\mathbb{A}_L}
\newcommand{\AAn}{\mathbb{A}_L\big|_{\Hom^\n}}
\newcommand{\B}{\mathcal{B}}
\newcommand{\Bk}{\textbf{B}_L^{\k}}
\newcommand{\Bnk}{\textbf{B}_L^{\k}\big|_{\Hom^\n}}
\newcommand{\Hom}{\mathcal{H}}
\newcommand{\NH}[2]{\Norme{#1}{\Sob^{#2}}}
\newcommand{\NP}[2]{\Norme{#1}{2,#2}}
\newcommand{\NHP}[2]{\NNorme{#1}{#2}}
\newcommand{\NR}[1]{\norme{#1}{}}

\begin{document}

\title{Normalizations with exponentially small remainders for
nonautonomous analytic periodic vector fields} 
\author{Tiphaine J\'ez\'equel \thanks{Institut de Math\'ematiques de Toulouse, Universit\'e Paul Sabatier, 118 route de Narbonne, 31062 Toulouse cedex 9, France.
Tel : 33 (0)5 61 55 65 87. 
E-mail : {\tt tiphaine.jezequel@math.univ-toulouse.fr}}}

\date{\today}
\maketitle
\begin{center}

\textit{Running head}: Normalizations with exponentially small remainders for periodic vector fields.

\begin{abstract}
In this paper we deal with  analytic nonautonomous vector fields with a periodic time-dependancy, that we study near an equilibrium point. In a first part,  we assume that the linearized system is split in two invariant subspaces $\Ezero$ and $\Eun$. Under light diophantine conditions on the eigenvalues of the linear part, we prove that there is a polynomial change of coordinates in $\Eun$ allowing to eliminate up to a finite polynomial order all terms depending only on the coordinate $\uzero\in\Ezero$ in the $\Eun$ component of the vector field. We moreover show that, optimizing the choice of the degree of the polynomial change of coordinates, we get an exponentially small remainder. In the second part, we prove  a normal form theorem with  exponentially small remainder. Similar theorems have been proved before in the autonomous case: this paper generalizes those results to the nonautonomous periodic case.

\vspace{1ex}

\end{abstract}
\end{center}

\section{Introduction} 
Let us consider an analytic differential system, near an equilibrium point that we take at the origin. To study the behaviour of solutions in this neighborhood, one can try to simplify the system by using a change of variables. "Simplify" can have several meanings: in fact what is expected is to obtain a system that we understand better than the initial system. Here we consider two different points of view. Precisely, let us consider an initial nonautonomous system of the form
\begin{equation}\label{syst}
\frac{du}{dt}=Lu+\V(u,t) , \quad t\in\R, \quad u(t)\in\R^\m,
\end{equation}
where $L$ is linear, and $\V$ is analytic, quadratic in $u$ and $T$-periodic, namely
\begin{equation}
\V(0,.)=0, \quad D_u\V(0,.)=0,\label{nonlin}
\end{equation}
\begin{equation}
\V(.,t+T)=\V(.,t)\text{ for all }t\in\R.\label{periodic}
\end{equation}
This means that we consider a nonautonomous system, but with an autonomous linear part. We will develop in this paper the two following ideas to simplify a such system.

\vspace{1ex}

\textbf{Simplification 1: Uncoupling subsets of coordinates; link with invariant manifolds.}
If we assume that our initial system can be split into 
\begin{equation}\label{introsyst}
\left\{
\begin{array}{rcl}
\Frac{d\uzero}{dt} &=& L_0\uzero+\V_0(\uzero,\uun,t), \\
  &&\\
  \Frac{d\uun}{dt} &=& L_1\uun+\V_1(\uzero,\uun,t),
\end{array}
\right.
\end{equation}

where $u=(\uzero,\uun)$, then one could want to find a change of variables of the form
$$\uun=\vun+\Phi(\uzero,t),$$
for which the system (\ref{introsyst}) is transformed into a new system 

\vspace{2ex}

\hfill
$
  \left\{\begin{array}{lr}
\Frac{d\uzero}{dt} =
L_0\uzero+\V^0(\uzero,\vun,t),& \hspace{12ex}\refstepcounter{equation}
(\theequation) \label{introeq1}\\
  \nonumber\\
  \Frac{d\vun}{dt} =
L_1\vun+\vun\V^1(\uzero,\vun,t).& \refstepcounter{equation}
(\theequation)\label{introeq2}
\end{array}\right.
$

\vspace{2ex}

 This idea of simplification is very close to the one used in KAM theory, for
Hamiltonian systems (see Kolmogorov (1954), Arnold (1978)). Here, when such a change of
variables exists, the main consequence is that $\vun(t)=0$ is a solution for
equation (\ref{introeq2}), and hence one can consider the reduced system
defined by equation (\ref{introeq1}). In particular, the set $\{\vun=0\}$ is an
invariant manifold for our differential system. In the case of an autonomous
system, this invariant manifold $\{\vun=0\}$ reads in the initial coordinates
$\{\uun=\Phi(\uzero)\}$. Here, working with periodic in time functions, this
manifold will be a periodic manifold $\{\uun=\Phi(\uzero,t)\}$. The search of
invariant manifolds and reduced systems is a key tool (Kelley (1967), for instance, develops this subject), widely used in the study of physical systems. For example, Haragus and Iooss (2010) provide lots of applications of the center manifold theorem.

\vspace{1ex}

\textbf{Simplification 2: Normal forms.} The normal form theories usually concern autonomous vector fields. In the case of nonautonomous analytic systems the same philosophy apply, if we consider the expansion of $\V$ in power series with time-dependant coefficients; the aim is then to find a change of variables of the form
$$u=\y+\Phi(\y,t),$$
such that our system (\ref{introsyst}) is transformed into a new system
\begin{equation}
\Frac{d\y}{dt}=L\y+\Nrond(\y,t),\nonumber
\end{equation}
in which $\Nrond$ is "as simple as possible". 
The original purpose of Poincaré in normal form theory was to obtain
$\Nrond=0$, but in general, when trying to eliminate monomials in the expansion
of $\V$, one can see that some resonant monomials remain whatever you do. So
$\Nrond$ will be an analytic function whose expansion in power series is only
made of the resonant monomials, which happen to be those satisfying a "normal
form criteria". There exist several normal form theories, leading to different
normal form criterias.  Here we work with the characterization introduced by Elphick \textit{et al.} (1987), i.e. we get power series made of monomials commuting with the
$\E^{tL^*}$ for all $t$ in $\R$. Precisely, we want $\Nrond$ to satisfy
\begin{equation}\label{critere}
\E^{-tL^*}\Nrond(\E^{tL^*}\y,t)=\Nrond(\y,0),\quad\text{ for all $\y\in\R^\m$ and 
all $t\in\R$}.
\end{equation}
Applying a normal form theorem, one expect that a system with less numerous 
monomials will be easier to study than the initial system which happens to be
often the case (see for instance Iooss and Adelmeyer (1992), Iooss and Peroueme (1993) and Lombardi (2000)).

\vspace{1ex}

\textbf{In fact, for both simplifications 1 and 2, it is very rarely possible to find exactly such changes of variables, but it is possible up to finite order.} Precisely, for a fixed integer $\p$, after changes of variables $\Phi_\p$ of degree $\p$, it was proved (see Haragus and Iooss (2010) for simplification 1, Iooss and Adelmeyer (1992) for simplification 2 in the autonomous case) that one can obtain transformed systems of the form
$$
\left\{
\begin{array}{rcl}
  \Frac{d\uzero}{dt} &=& L_0\uzero+\V^0(\uzero,\vun,t),\\
  &&\\
  \Frac{d\vun}{dt} &=& L_1\vun+\vun\V^1_\p(\uzero,\vun,t)+\Reste_\p(\uzero,t),
\end{array}
\right.$$
for simplification 1, and
\begin{equation}
\frac{d\y}{dt}=L\y+\Nrond_\p(\y,t)+\Reste_\p(\y,t),\nonumber
\end{equation}
for simplification 2; where $\V^1_\p$ and $\Nrond_\p$ are polynomials of degree 
$\p$, and $\Reste_\p$ is analytic of order larger than $\p$. It is an
interesting result since a usual way to study the dynamics of a vector field is
to study the system truncated at a fixed order in the expansion in power
series, and then to consider the complete vector field as a perturbation of
this truncated system (see Guckenheimer and Holmes (1983),Iooss and Peroueme (1993), Lombardi (2000), Chossat and Iooss (1994)).

But to apply perturbation theory, those results are particularly interesting 
if the remaining part $\Reste_\p$ happens to be small. So the next key idea is
to optimize the degree $\p$ of the truncation, to minimize the size of the
remaining part $\Reste_\p$. For autonomous systems, Iooss and Lombardi (2010, for
transformation 1 and 2005, for transformation 2) followed
an idea developped before for Hamiltonian systems (developped by Nekoroshev (1977, 1979); see also Pöschel (1993)) to prove that when some light hypothesis
hold, the remaining part $\Reste_{\p_{opt}}$ can be found exponentially small. 

\vspace{1ex}

\textbf{In this paper, we generalize those theorems with exponentially small 
remainders to nonautonomous systems with a periodic time dependancy and an
autonomous linear part.} We prove here in our theorem \ref{th1} (for
simplification 1) and our theorem \ref{thm2} (for simplification 2) that, whith
light hypothesis of non-resonance on the non-linear part, we have the same kind
of results that for autonomous systems; in particular, notice that we obtain
time-independant exponentially small bounds for the remainder.

\vspace{1ex}

\textbf{This situation arises quite often in applications}, when a system is 
periodically forced. Our theorem \ref{th1} and its corollary \ref{cor1} might
be used for instance in the case considered by Touzé
and Amabili (2006) in which they build reduced-order models for
nonlinear vibrations of structures. In their section 3.2, they consider a
two-dof system with damping terms ($\xi_1$ and $\xi_2$) and with a periodic
forcing ($F_1$), of the form

\vspace{2ex}

\hfill
$
\left\{\begin{array}{lr}
  \frac{d^2 X_1}{dt^2}+2\xi_1\omega_1\frac{dX_1}{dt}+\omega_1^2
X_1=V_1(X_1,X_2)+F_1 \cos(\Omega t),&\hspace{10ex}\refstepcounter{equation}
(\theequation)\label{dof1}\\
  &\nonumber\\ 
  \Frac{d^2 X_2}{dt^2}+2\xi_2\omega_2\Frac{dX_2}{dt}+\omega_2^2
X_2=V_2(X_1,X_2)&\hspace{10ex}\refstepcounter{equation}
(\theequation)\label{dof2}
\end{array}\right.
$

\vspace{2ex}

where $V_1,V_2$ are quadratic in $(X_1,X_2)$. In their analysis, they consider that $\xi_1$ is small and
then that  the first oscillator $X_1$ plays the role of a central manifold:
namely they use the central manifold theorem to state that the mode $X_2$ does
not awake as $t$ goes to infinity, thus they make the approximation $X_2=0$ to
solve the first equation (\ref{dof1}) and then solve the second equation
(\ref{dof2}) with the value $X_1(t)$ computed. They observe that this process
gives good qualitative results for small values of $F_1$. In fact, the center
manifold theorem does not apply here because this system is not autonomous when
$F_1$ is nonzero. Here we show in our corollary \ref{cor1} that if this forcing term $F_1$ is small, then the invariant manifold nearly remains time-independant. Precisely, if one take $F_1=\eps^2$, then the corollary \ref{cor1} ensures that the equation leading to the invariant manifold reads
$$X_2=\Phi_A(X_1)+\eps \Phi_B(X_1,\eps,t)$$
where $\Phi_A$ does not depend on $t$.

\vspace{2ex}

In section \ref{main} we gather the main results of the paper. We need for that purpose to introduce, in subsection \ref{notations}, a few notations. We then state in subsection \ref{subth1} our first main theorem (theorem \ref{th1}), which deals with the "Simplification 1" described above, and we state in subsection \ref{subth2} the second main theorem (theorem \ref{thm2}) which is a normal form theorem ("Simplification 2" above). Then the rest of the paper is devoted to the proofs of these two theorems: we detail the proof of theorem \ref{th1} in section \ref{proof1} (the strategy of proof is introduced in the first subsection (\ref{sketch})), and in section \ref{proof2}, we give the main ideas of proof of theorem \ref{thm2} (here also the strategy is given in the first subsection \ref{constr}).

\section{Notations and main results} \label{main}

We gather in this section the main theorems proved in this paper. In the whole paper, we consider a differential system of the form (\ref{syst}), assuming that (\ref{nonlin}) holds. We also assume that the map $\V$ is analytic in $u$, and that its time-dependancy is $\T$-periodic with regularity $\Sob^\l$. For that purpose we need to define precisely a space $\mathcal{A}(\R^\m,\Sob^\l(\R/\T\Z, \R^\m))$ of such functions. So, in the following subsection, we begin by defining this set of functions and other usefull sets.

\subsection{A few notations}\label{notations}
We first recall what we denote by $\Sob^\l(\R/\T\Z, \R^\m)$.
\begin{defi}
$\Sob^\l(\R/\T\Z, \R^\m)$ stands for the Sobolev space of functions $f$ from $\R/\T\Z$ to $\R^\m$ whose Fourier coefficients $f^{(\k)}\in\R^\m$ satisfy
$$\sum_{\k=-\infty}^{+\infty}(1+\k^2)^\l\Norme{f^{(\k)}}{\R^\m}^2<+\infty.$$
\end{defi}
We can then define the space we wanted for $\V$.
\begin{defi}
For any neighborhood $\Nu$ of the origin in $\R^\m$ we denote by $\mathcal{A}(\Nu,\Sob^\l(\R/\T\Z, \R^\m))$ the space of maps $\V$ for which there exists a family of $q$-linear symmetric maps $(\V_q(t))_{q\geq0}$ on $(\R^m)^q$, with a radius of convergence $\rho$ and a positive constant $\c$ such that
$$\V(u,t)=\sum^{+\infty}_{q=0}\V_q(t)[u^{(q)}],$$
(here $[u^{(q)}]$ stands for the $q$-uple of vectors $[u,u,\cdots,u]$) and 
\begin{equation} 
\NH{\V_q(.)[u_1,\cdots,u_q]}{\l}\leq \frac{\c}{\rh^q}|u_1||u_2|\cdots|u_q|\quad\text{for all }u_1,\cdots,u_q\in\R^\m.\nonumber
\end{equation}
\end{defi}
We also need, to state our main theorems, the following spaces.
\begin{defi}
Let $\mathcal{P}_{\p}(\R^{\m_0}, \Sob^{\l}(\R/\T\Z,\R^{\m_1}))$ be the space of polynomials $P$ of degree less than $\p$, namely which read
$$P(x_1,\cdots,x_{\m_0},t)=\sum_{|\alpha|\leq\p}P_\alpha(t)x_1^{\alpha_1}\cdots x_{\m_0}^{\alpha_{\m_0}}$$
where $\alpha=(\alpha_1,\cdots,\alpha_{\m_0})$ belongs to $\N^{\m_0}$, $P_\alpha$ to $\Sob^\l(\R/\T\Z, \R^{\m_1})$, and
$$|\alpha|=|\alpha_1|+\cdots+|\alpha_{\m_0}|.$$
And let $\Hom_{\n}(\R^{\m_0}, \Sob^{\l}(\R/\T\Z,\R^{\m_1}))$ be the space of homogeneous polynomials $P$ of degree $\n$, namely of the form
$$P(x_1,\cdots,x_{\m_0},t)=\sum_{|\alpha|=\n}P_\alpha(t)x_1^{\alpha_1}\cdots x_{\m_0}^{\alpha_{\m_0}}$$
where the $P_\alpha$ are in $\Sob^\l(\R/\T\Z, \R^{\m_1})$.
\end{defi}

\subsection{Uncoupling subsets of coordinates}\label{subth1}
In this first subsection, we state a theorem where a change of variables uncouple a subset of variables from another one (it is the "Simplification 1" of the introduction). We suppose here that the linear operator $L$ is the direct sum of two linear operators. Precisely, we assume
\begin{hyp}\label{hyp2}
Assume that 
\begin{enumerate}
\item there exist a neighborhood $\Nu$ of the origin in $\R^\m$ and an integer $\l\geq1$ such that $\V$ belongs to $\mathcal{A}(\Nu,\Sob^\l(\R/\T\Z, \R^\m))$;
\item $L$ is the direct sum of the linear operators $L_0$ on $\Ezero$ (dimension $\m_0$) and $L_1$ on $\Eun$ (dimension $\m_1$), with  $\Ezero\oplus\Eun=\R^\m$ and $L_0$ diagonizable. Hence, spliting $u$ in $u=\uzero+\uun$  and $\V=\V_0+\V_1$ on $\Ezero\oplus\Eun$, the system (\ref{syst}) reads
\begin{equation}\label{systsplit}
\left\{
\begin{array}{rcl}
\Frac{d\uzero}{dt} &=& L_0\uzero+\V_0(\uzero,\uun,t), \\
  &&\\
  \Frac{d\uun}{dt} &=& L_1\uun+\V_1(\uzero,\uun,t),
\end{array}
\right.
\end{equation}
\item there exist two positive constants $\ga$ and $\ta$ such that, denoting the eigenvalues of $L_0$ by $\lambda_1^{(0)},\cdots,\lambda_{\m_0}^{(0)}$ and by $\lambda_1^{(1)},\cdots,\lambda_{\m_1}^{(1)}$ the eigenvalues of $L_1$, then for all $\a$ in $\mathbb{N}^{\m_0}$, all $\k$ in $\mathbb{Z}$ and all $\j$, $1\leq\j\leq\m_1$,
\begin{equation} \label{diophantian}
|\scal{\a}{\lambda^{(0)}}+\I\k \frac{2\pi}{\T}-\lambda_j^{(1)}|\geq \frac{\ga}{(|\a|+|k|)^\ta}
\end{equation}
holds. 
\item let $\nu$ be the maximal size of the Jordan blocks of $L_1$, then
\begin{equation} \label{l}
\ta\nu\leq\l.
\end{equation}
\end{enumerate}
\end{hyp}
\begin{rem}
We need to take $\l\geq1$ to insure that $\Sob^\l$ is an algebra.
\end{rem}
\begin{rem}
The hypothesis (\ref{diophantian}) is called a "non-resonance" hypothesis: without taking $\k$ into account, this hypothesis means that the eigenvalues of $L_1$ cannot happen to be a sum of eigenvalues of $L_0$, and that moreover we have an estimate of how different they are. Taking $\k$ into account, it expresses that, even after periodic translations of time, non-resonance remains.
\end{rem}
We can now state
\begin{thm} \label{th1}
Consider the system (\ref{systsplit}), and suppose that hypothesis \ref{hyp2} holds.
Then there exists $\delta_0>0$ such that, for all $\delta$ in $]0,\delta_0[$ there exist an integer $\p_{\delta}$ with 
\begin{equation}
\p_{\delta}=\grandO{\delta\rightarrow 0}(\delta^{-\b}),\nonumber
\end{equation}
and a function $\Phi_\delta$ in $\mathcal{P}_{\p_\delta}(\Ezero, \Sob^{\l+1}(\R/\T\Z,\Eun))$ (where $\mathcal{P}_{\p}$ was defined in section \ref{notations}) with
$$\Phi_\delta(0,.)=0, \quad D_{\uzero}\Phi_\delta(0,.)=0;$$
such that the change of variables
\begin{equation} \label{cdv}
\uun=\vun+\Phi_\delta(\uzero,t)
\end{equation}
 transforms the system (\ref{systsplit}) into

\begin{equation}\label{systnorm}
\left\{
\begin{array}{rcl}
\Frac{d\uzero}{dt} &=& L_0\uzero+\V^0(\uzero,\vun,t),\\
  &&\\
\Frac{d\vun}{dt}&=& L_1\vun+\V^1(\uzero,\vun,t)+\Reste(\uzero,t);
\end{array}
\right.
\end{equation}

in which $\V^0,\V^1$ and $\Reste$ are analytic, and satisfy
\begin{equation}\label{vun}
\NH{\V^1(\uzero,\vun,.)}{\l}\leq\Mzero\NR{\vun}(\NR{\uzero}+\NR{\vun}) \quad \text{for }\NR{\uzero},\NR{\vun}\leq\delta_0,
\end{equation}
and
\begin{equation}\label{maj}
\supl_{\NR{\uzero}\leq\delta}\NH{\Reste(\uzero,.)}{\l}\leq\M\E^{-\frac{\om}{\delta^\b}},
\end{equation}
where $\M,\Mzero,\om$ depend only on $\T, \m_0, \m, \c, \rh, L, \l,\delta_0$ and $\ta$, and where
$$\b=\TFrac{1}{\l+\ta\nu+1}.$$
Moreover, $\V^0$ reads $\V^0(\uzero,\vun,t)=\V_0(\uzero+\vun+\Phi_\delta(\uzero,t),t).$
\end{thm}
\begin{rem}
Observe that with the system in the new form (\ref{systnorm}) $\vun=0$ is "very close" to solve the second equation, since (\ref{vun}) guarantees that for $\vun=0$ we have $\V^1(\uzero,0,.)=0$, and since (\ref{maj}) ensures that $\Reste(\uzero,.)$ remains exponentially small. Then this theorem expresses that the manifold $\{\uun=\Phi_\delta(\uzero,t)\}$ is "exponentially close" to be an invariant manifold for our system (\ref{syst}). 
\end{rem}
\begin{rem}
This theorem is the generalization to periodic time-dependant vector fields of the Theorem 1 of Iooss and Lombardi (2010). Notice that, unlike the latter, here $\M$ and $\om$ also depend on the dimension $\m$, and not only on the dimension $\m_0$ of $\Ezero$, so that we cannot consider systems of infinite dimension.
\end{rem}

Moreover, we have the following proposition, which gives more precision about what happens if the time-dependancy appear as a small perturbation of an autonomous system.

\begin{cor}\label{cor1}
Consider the system (\ref{systsplit}) and suppose that hypothesis \ref{hyp2} holds. Assume that $\V$ also depends analytically on a parameter $\eps$ in the following way: 
\begin{equation}
\V_\eps(u,t)=A(u)+\eps B(u,\eps,t)+\eps^2 C(\eps,t)\label{Veps}
\end{equation}
with $B(0,.,.)=0$. 

Setting $U:=(u,\eps)$, if $A,B$ and $C$ belong to $\mathcal{A}(\R^{\m+1},\Sob^\l(\R/\T\Z, \R^\m))$ then the theorem \ref{th1} apply and moreover
$$\Phi(\uzero,\eps,t)=\Phi_A(\uzero)+\eps\Phi_{B,C}(\uzero,\eps, t)$$
where $\Phi_A$ is the change of variables computed applying theorem \ref{th1} to our system (\ref{systsplit}) at $\eps=0$.
\end{cor}

\begin{rem}
Remark that $(\ref{Veps})$ prevents us to take a periodic forcing of the form 
$A(u)+\eps F(t)$. We need  hypothesis $(\ref{Veps})$, to apply theorem
\ref{th1} because we want hypothesis $D_U V(U=0,.)=0$ to hold (see
(\ref{nonlin})).

But we can take a periodic forcing $A(u)+\eps^2 F(t)$, and then obtain that the periodic invariant manifold $\{\uun=\Phi_A(\uzero)+\eps\Phi_{B,C}(\uzero,\eps, t)\}$ is $\eps$-close to the autonomous invariant manifold $\{\uun=\Phi_A(\uzero)\}$.
\end{rem}

\begin{rem}
The proof of this corollary directly follows from the proof of theorem \ref{th1}, the details are left to the reader.
\end{rem}

\subsection{Normal form}\label{subth2}
In this second subsection, we state a normal form theorem (it is the "Simplification 2"' in the introduction). We only assume that the following ''non-resonance'' hypothesis holds:
\begin{hyp}\label{hyp3}
Assume that
\begin{enumerate}
\item there exist a neighborhood $\Nu$ of the origin in $\R^\m$ and an integer $\l\geq1$ such that $\V$ belongs to $\mathcal{A}(\Nu,\Sob^\l(\R/\T\Z, \R^\m))$;
\item $L$ is diagonizable;
\item there exist two positive constants $\ga$ and $\ta$ such that, denoting the eigenvalues of $L$ by $\lambda_1,\cdots,\lambda_\m$ in $\CC^\m$, then for all $\a$ in $\mathbb{N}^{\m_0}$, all $\k$ in $\mathbb{Z}$ and all $\j$, $1\leq\j\leq\m$
\begin{equation} \label{diophantian2}
|\scal{\a}{\lambda}+\I\k \frac{2\pi}{\T}-\lambda_j|\geq \frac{\ga}{(|\a|+|k|)^\ta};
\end{equation}
\item $\ta\leq\l.$
\end{enumerate}
\end{hyp}

We then have the following normal form theorem.
\begin{thm}\label{thm2}
Consider the system (\ref{syst}), and suppose that hypothesis \ref{hyp3} 
holds. Then for all $\delta>0$ there exist an integer $\p_{\delta}$ with 
\begin{equation}
\p_{\delta}=\grandO{\delta\rightarrow 0}(\delta^{-b}),\nonumber
\end{equation}
and a function $\Phi_\delta$ in $\mathcal{P}_{\p_\delta}(\R^\m, \Sob^{\l+1}(\R/\T\Z,\R^\m))$ (where $\mathcal{P}_{\p}$ is defined in section \ref{notations}) with 
$$\Phi_\delta(0,.)=0, \quad D_u\Phi_\delta(0,.)=0,$$
such that the change of variables
\begin{equation}\label{NFcdv}
u=\y+\Phi_\delta(\y,t)
\end{equation}
transforms the system (\ref{syst}) into the normal form
\begin{equation}\label{systNF}
\frac{d\y}{dt}=L\y+\Nrond(\y,t)+\Reste(\y,t)
\end{equation}
where $\Nrond$ belongs to $\mathcal{P}_{\p_\delta}(\R^\m, \Sob^{\l}(\R/\T\Z,\R^\m))$ and satisfy
$$\Nrond(0,.)=0, \quad D_u\Nrond(0,.)=0,$$
and the normal form criteria 
\begin{equation}\label{critereNF}
\E^{-tL^*}\Nrond(\E^{tL^*}\y,t)=\Nrond(\y,0),\quad\text{ for all $\y\in\R^\m$ and all $t\in\R$};
\end{equation}
and where the remainder $\Reste$ is analytic and satisfies
\begin{equation}\label{NFmaj}
\supl_{\NR{\y}\leq\delta}\NH{\Reste(\y,.)}{\l}\leq\M'\delta^2\E^{-\frac{\om}{\delta^{\b}}},
\end{equation}
with $\M'$ and $\om$ depending only on $\T, \m, \c, \rh, L, \l,\delta_0$ and $\ta$, and 
$$\b=\TFrac{1}{\l+\ta+1}.$$
\end{thm}
\begin{rem}
This theorem is a typical normal form theorem: its aim is to simplify the initial system with the aid of a change of variables. Indeed, the polynomial change of variables (\ref{NFcdv}) transforms the system (\ref{syst}) into a new system (\ref{systNF}) in which the polynomial part $\Nrond$ is simpler: the normal form criteria (\ref{critereNF}) means that in $\Nrond$, all the monomials which does not commute with all the $\E^{sL^*}$ have been eliminated by the change of variable.

Moreover, while for some given $\p$ the remainder is polynomially small, this theorem optimizes the degree $\p=\p_{opt}(\delta)$ of the polynomial part $\Nrond$, so that the system is nearly reduced to a polynomial system, since (\ref{NFmaj}) ensures that the remainder $\Reste$ is exponentially small.
\end{rem}
\begin{rem}
This theorem is the generalization to periodic time-dependant vector fields of the Theorem 1.4 of Iooss and Lombardi (2005).
\end{rem}

\section{Proof of theorem \ref{th1}}\label{proof1}

This section is entirely devoted to the proof of theorem \ref{th1}.

\subsection{Strategy of proof}\label{sketch}
First, fix $\delta>0$ and $\p$ in $\N$, and see later what conditions on $\delta$ and $\p$ need to be satisfied. 

One can check that, for any function $\Phi$, the change of variables (\ref{cdv})
\begin{equation}
\uun=\vun+\Phi(\uzero,t)\nonumber
\end{equation}
transforms the system (\ref{systsplit}) into the new system (\ref{systnorm}) with
\begin{eqnarray}
\V^0(\uzero+\vun,t)&=&\V_0(\uzero+\vun+\Phi(\uzero,t),t)\label{eq1}, \\
\V^1(\uzero+\vun,t)+\Reste(\uzero,t)&=& \V_1(\uzero+\vun+\Phi(\uzero,t),t)\label{eq2}\\
																		&&-D_{\uzero}\Phi(\uzero,t).\V_0(\uzero+\vun+\Phi(\uzero,t),t)\nonumber\\
																		&&-(D_{\uzero}\Phi(\uzero,t).L_0\uzero-L_1\Phi(\uzero,t))\nonumber\\
																		&&- \partial_t\Phi(\uzero,t).\nonumber
\end{eqnarray}
Thus, for a fixed $\Phi$, equation (\ref{eq1}) provides the value of $\V^0$. Then, we look for a function $\Phi$ such that $\V^1$ and $\Reste$, whose sum is computed in (\ref{eq2}), satisfy
\begin{eqnarray}
\NH{\V^1(\uzero+\vun,.)}{\l} \leq & \Mzero\NR{\vun}(\NR{\uzero}+\NR{\vun}), &\text{ when } \NR{\uzero}\leq\delta_0, \label{majvun}\\
\NH{\Reste(\uzero,.)}{\l} \leq & \M\E^{-\frac{\om}{\delta^\b}}, &\text{ when } \NR{\uzero}\leq\delta \label{majexp}.
\end{eqnarray}
Necessarily, if (\ref{majvun}) holds, then for $\vun=0$ equation (\ref{eq2}) becomes
\begin{eqnarray}
\Reste(\uzero,t)&=& \V_1(\uzero+\Phi(\uzero,t),t)-D_{\uzero}\Phi(\uzero,t).\V_0(\uzero+\Phi(\uzero,t),t)\nonumber\\
								&&-(D_{\uzero}\Phi(\uzero,t).L_0\uzero-L_1\Phi(\uzero,t))- \partial_t\Phi(\uzero,t). \label{eqprime}
\end{eqnarray}
Define
\begin{equation}
(\A_L\Phi)(\uzero,t):=D_{\uzero}\Phi(\uzero,t).L_0\uzero-L_1\Phi(\uzero,t).\nonumber
\end{equation}
Then (\ref{eqprime}) reads
\begin{eqnarray}
(\A_L+\partial_t)\Phi(\uzero,t)+\Reste(\uzero,t)&=&\V_1(\uzero+\Phi(\uzero,t),t)\nonumber\\
																&&-D_{\uzero}\Phi(\uzero,t).\V_0(\uzero+\Phi(\uzero,t),t).\nonumber
\end{eqnarray}
Let us denote by $\Pi_\p$ the projection of $\mathcal{A}(\Ezero,\Sob^{\l}(\R/\T\Z,\Eun))$ on 

$\mathcal{P}_{\p}(\Ezero,\Sob^{\l}(\R/\T\Z,\Eun))$. 

\vspace{2ex}

\textbf{Here is the strategy of proof, in three steps, that we will follow in the next subsections}:
\begin{description}
\item[Step A]: In section \ref{construct}, for any fixed $\p$, we prove the existence of some $\Phi$ in $\mathcal{P}_{\p}(\Ezero,\Sob^{\l+1}(\R/\T\Z,\Eun))$ such that
\begin{eqnarray}
(\A_L+\partial_t)\Phi(\uzero,t)&=&\Pi_{\p}(\V_1(\uzero+\Phi(\uzero,t),t)\label{step1}\\
																&&-D_{\uzero}\Phi(\uzero,t).\V_0(\uzero+\Phi(\uzero,t),t)).\nonumber
\end{eqnarray}
\item[Step B]: In section \ref{upperbound}, using the $\Phi$ computed in Step A, we set
\begin{eqnarray}
\Reste(\uzero,t)&=&(\Id-\Pi_{\p})(\V_1(\uzero+\Phi(\uzero,t),t)\label{step2}\\
										&&-D_{\uzero}\Phi(\uzero,t).\V_0(\uzero+\Phi(\uzero,t),t)).\nonumber
\end{eqnarray}
Thus, $\Reste(\uzero,t)=\grandO{}(\NR{\uzero}^{\p+1})$. We then compute upper bounds of $\Reste$ of the form
\begin{equation}
\NH{\Reste(\uzero,.)}{\l}\leq\M_{\p}\NR{\uzero}^{\p+1},\nonumber
\end{equation}
with Gevrey estimates for $\M_{\p}$. And then, for a given $\delta$, we chose an integer $\p=\p_{opt}$ minimizing $\M_{\p}\delta^{\p+1}$. We will see that if $\delta$ is sufficiently small, for our choice $\p_{opt}(\delta)$, the estimate (\ref{majexp}) holds.
\item[Step C]: From (\ref{eq2}), we get the value of $\V^1$: 
\begin{eqnarray}
\V^1(\uzero,\vun,t)&:=&-\Reste(\uzero,t)+\V_1(\uzero+\vun+\Phi(\uzero,t),t)\nonumber\\
										&&-D_{\uzero}\Phi(\uzero,t).\V_0(\uzero+\vun+\Phi(\uzero,t),t)\nonumber\\
										&&-(\mathcal{A}_L+\partial_t)\Phi(\uzero,t).\nonumber
\end{eqnarray}
And it remains to show, in section \ref{StepC}, that there exists $\Mzero$ such that (\ref{majvun}) holds.
\end{description}
But first, to facilitate the estimates in Step B, we introduce apropriate norms on $\mathcal{P}_{\p}(\Ezero,\Sob^{\l}(\R/\T\Z,\Eun))$.

\subsection{Norms on $\mathcal{P}_{\p}(\Ezero,\Sob^{\j}(\R/\T\Z,\R^{\m}))$}\label{defnorm}
In fact, we define norms on the spaces of homogeneous polynomials of degree $\n$ $\Hom^{\n}(\Ezero,\Sob^{\j}(\R/\T\Z,\R^{\m}))$ (defined in section \ref{notations}). Let $(e_1,\cdots,e_{\m_0})$ be a basis of $\Ezero$ and $(e_{\m_0+1},\cdots,e_{\m})$ a basis of $\Eun$. Then, if $f$ belongs to $\Hom^{\n}(\Ezero,\Sob^{\j}(\R/\T\Z,\R^{\m}))$, it reads
\begin{equation}
f(\uzero,t)=\suml_{\dindice{\NR{\alpha}=\n}{i=1,\cdots,\m}} f_{\alpha,i}(t)x_1^{\alpha_1}\cdots x_{\m_0}^{\alpha_{\m_0}}e_i:=\suml_{\NR{\alpha}=\n} f_{\alpha}(t)x_1^{\alpha_1}\cdots x_{\m_0}^{\alpha_{\m_0}},\nonumber
\end{equation}
where 
$$\uzero:=x_1e_1+\cdots+x_{\m_0}e_{\m_0};$$
and $f_{\alpha,i}\in\Sob^{\j}(\R/\T\Z,\R))$, $f_{\alpha}\in\Sob^{\j}(\R/\T\Z,\R^{\m})$.

\begin{defi}
If $f$ belongs to $\Hom^{\n}(\Ezero,\Sob^{\j}(\R/\T\Z,\R^{\m}))$ and reads
$$f(\X,t)=\suml_{\NR{\alpha}=\n}f_{\alpha}(t)\X^{\alpha},$$ 
\vspace{-1ex}
where $$\X:=(\X_1,\cdots,\X_{\m_0}),\quad \X^{\alpha}:=\X_1^{\alpha_1}\cdots\X_{\m_0}^{\alpha_{\m_0}};$$ then define
\begin{equation}
\NHP{f}{\n,\Sob^\j}^2\hspace{-1ex}:=\NP{\suml_{\NR{\alpha}=\n}\NH{f_{\alpha}}{\j}\X^{\alpha}}{\n}^2\hspace{-2ex}=\suml_{\NR{\alpha}=\n}\frac{\alpha!}{n!}\NH{f_{\alpha}}{\j}^2=\hspace{-1ex}\sum^{+\infty}_{\k=-\infty}(1+\k^2)^{\j}\NP{\suml_{\NR{\alpha}=\n}f_{\alpha}^{(\k)}\X^{\alpha}}{\n}^2\hspace{-2ex};\nonumber
\end{equation}
where $\alpha!:=\alpha_1!\cdots\alpha_{\m_0}!$, and
\begin{itemize}
\item $\NP{.}{\n}$ is the norm on $\Hom^\n(\R^{\m_0},\R^{\m_1})$, introduced by Iooss and Lombardi (2005), defined by
$$\NP{P(\X)}{\n}^2=\NP{\suml_{\NR{\alpha}=\n}P_{\alpha}\X^{\alpha}}{\n}^2:=\suml_{\NR{\alpha}=\n}\frac{\alpha!}{n!}\Norme{P_{\alpha}}{\R^{\m_1}}^2;$$
\item $f_{\alpha}^{(\k)}\in\R^\m$ stands for the $\k^{\text{th}}$ Fourier coefficient of $f_\alpha$;
\item $\NH{.}{\j}$ is the canonical norm on $\Sob^{\j}(\R/\T\Z, \R^{\m})$
$$\NH{f}{\j}^2:=\sum^{+\infty}_{\k=-\infty}(1+\k^2)^{\j}\Norme{f_{\alpha}^{(\k)}}{\R^\m}^2.$$
\end{itemize}
In particular, to simplify notations, we denote
$$\NHP{f}{\n}:=\NHP{f}{\n,\Sob^{\l}}.$$
\end{defi}

\subsection{Step A : construction of $\Phi$}\label{construct}
In this subsection we fix an integer $\p$, and our aim is to construct $\Phi$ in $\mathcal{P}_{\p}(\Ezero,\Sob^{\j}(\R/\T\Z,\Eun))$ such that (\ref{step1}) is satisfied. Denote $$\Phi(\X,t):=\sum^{\p}_{\n=2}\Phi_{\n}(\X,t)$$ where $\Phi_{\n}$ belongs to $\Hom^{\n}(\Ezero,\Sob^{\l+1}(\R/\T\Z,\Eun))$. We project (\ref{step1}) on the spaces $\Hom^{\n}(\Ezero,\Sob^{\l}(\R/\T\Z,\Eun))$, for $\n\leq\p$. Denoting by $\pi_{\n}$ this projection on 

$\Hom^{\n}(\Ezero,\Sob^{\l}(\R/\T\Z,\Eun))$, we obtain:
\begin{eqnarray} 
(\A_L+\partial_t)\Phi_\n(\uzero,t)&=&\pi_{\n}(\V_1(\uzero+\Phi(\uzero,t),t)\label{projn}\\
&&\hspace{5ex}-D_{\uzero}\Phi(\uzero,t).\V_0(\uzero+\Phi(\uzero,t),t)).\nonumber
\end{eqnarray}
Expanding the right hand side of (\ref{projn}) in power series, one can observe that, since
$$\V(0,.)=0=D_{\uzero}\V(0,.),$$
(because of (\ref{nonlin})), this right hand side of (\ref{projn}) only depends on $\Phi_2,\cdots,\Phi_{\n-1}$ and $\uzero,t$. So, if $(\A_L+\partial_t)$ is invertible, then (\ref{projn}) enables us to construct the $\Phi_{\n}$ for $2\leq\n\leq\p$ by induction. The rest of this subsection is devoted to the proof of the invertibility of $\AAL:=(\A_L+\partial_t)$. Precisely, we prove the following

\begin{lem}\label{leminv}
If hypothesis \ref{hyp2} holds, then the linear operator
\begin{eqnarray}
\AAn:\Hom^{\n}(\Ezero,\Sob^{\l+1}(\R/\T\Z,\Eun))&\longrightarrow&\Hom^{\n}(\Ezero,\Sob^{\l}(\R/\T\Z,\Eun))\nonumber\\
																				\Phi &\longmapsto& (\A_L+\partial_t)\Phi \nonumber
\end{eqnarray}
is invertible.

Moreover, for all $\j$, $0\leq\j\leq\l+1$ and all $\F$ in $\Hom^{\n}(\Ezero,\Sob^{\l}(\R/\T\Z,\Eun))$ we have 
\begin{equation}\label{majinv}
\NHP{\AAn^{-1}(\F)}{\n,\Sob^{\j}}\leq \C_\j\n^{\j+\ta\nu}\NHP{\F}{\n}
\end{equation}
where $\nu$ is defined in hypothesis \ref{hyp2} and
\begin{eqnarray}
\Lambda&:=&\max\{\NR{\lambda_i^{(0)}},\NR{\lambda_j^{(1)}}, 1\leq i\leq \m_0, 1\leq j\leq\m_1\},\nonumber\\
\C_\j&:=&\max(1,\frac{\nu}{\ga^{\nu}})(1+\frac{\T^2}{2\pi^2}(4\Lambda^2+1))^{\frac{\j}{2}}.\nonumber
\end{eqnarray}

\end{lem}

\textbf{Proof.}
To prove this lemma, we use Fourier coefficients, so that instead of one operator in $\Hom^{\n}(\Ezero,\Sob^{\l}(\R/\T\Z,\Eun))$, we deal with an infinity of linear operators in $\Hom^{\n}(\Ezero,\Eun)$. Then, for these linear operators we can use the results of invertibility stated by Iooss and Lombardi (2010). 

Indeed, any function $\Phi$ in $\Hom^{\n}(\Ezero,\Sob^{\l+1}(\R/\T\Z,\Eun))$ reads
$$\Phi(\X,t):=\suml_{\NR{\alpha}=\n}\phi_{\alpha}(t)\X^{\alpha}$$
with $\phi_\alpha$ in $\Sob^{\l+1}(\R/\T\Z,\Eun)$. Since $\l+1\geq1$, we have the Fourier expansions
\begin{eqnarray}
\phi_\alpha(t)&=&\sum^{+\infty}_{\k=-\infty}\phi_{\alpha}^{(\k)}\E^{\I\k\frac{2\pi}{\T}t}\nonumber\\
\partial_t\phi_\alpha(t)&=&\sum^{+\infty}_{\k=-\infty}\I\k\frac{2\pi}{\T}\phi_{\alpha}^{(\k)}\E^{\I\k\frac{2\pi}{\T}t},\nonumber\\
\mathcal{A}_L\phi_\alpha(t)&=&\sum^{+\infty}_{\k=-\infty}\mathcal{A}_L\phi_{\alpha}^{(\k)}\E^{\I\k\frac{2\pi}{\T}t}.\nonumber
\end{eqnarray}
Thus, we obtain, introducing the notation $\Phi^{(\k)}(\X)$
\begin{eqnarray}
\Phi(\X,t)&=&\suml_{\NR{\alpha}=\n}\big(\sum^{+\infty}_{\k=-\infty}\phi_{\alpha}^{(\k)}\E^{\I\k\frac{2\pi}{\T}t}\big)\X^\alpha=\sum^{+\infty}_{\k=-\infty}\big(\suml_{\NR{\alpha}=\n}\phi_{\alpha}^{(\k)}\X^\alpha\big)\E^{\I\k\frac{2\pi}{\T}t}\nonumber\\
				&:=&\sum^{+\infty}_{\k=-\infty} \Phi^{(\k)}(\X)\E^{\I\k\frac{2\pi}{\T}t};\nonumber
\end{eqnarray}
and 
\begin{eqnarray}
\AAL\Phi(\X,t)&=&(\mathcal{A}_L+\partial_t)(\suml_{\NR{\alpha}=\n}\phi_{\alpha}(t)\X^{\alpha})\nonumber\\						&=&\suml_{\NR{\alpha}=\n}\left(\sum^{+\infty}_{\k=-\infty}(\mathcal{A}_L+\I\k\frac{2\pi}{\T})\phi_{\alpha}^{(\k)}\E^{\I\k\frac{2\pi}{\T}t}\right)\X^{\alpha}\nonumber\\					&=&\sum^{+\infty}_{\k=-\infty}\left(\suml_{\NR{\alpha}=\n}(\mathcal{A}_L+\I\k\frac{2\pi}{\T})\phi_{\alpha}^{(\k)}\X^\alpha\right)\E^{\I\k\frac{2\pi}{\T}t}\nonumber\\
&=&\sum^{+\infty}_{\k=-\infty}\left((\mathcal{A}_L+\I\k\frac{2\pi}{\T})\suml_{\NR{\alpha}=\n}\phi_{\alpha}^{(\k)}\X^\alpha\right)\E^{\I\k\frac{2\pi}{\T}t}\nonumber\\
&=&\sum^{+\infty}_{\k=-\infty}(\mathcal{A}_L+\I\k\frac{2\pi}{\T})\Phi^{(\k)}(\X)\E^{\I\k\frac{2\pi}{\T}t}.\nonumber					
\end{eqnarray}
We then proceed in several steps.

\vspace{2ex}

\textbf{Step 1: We first show that the operator $\AAn$ is injective.}

We fix $F$ in $\Hom^{\n}(\Ezero,\Sob^{\l}(\R/\T\Z,\Eun))$. Denote
\begin{eqnarray}
F(\X,t)&:=&\suml_{\NR{\alpha}=\n}F_{\alpha}(t)\X^{\alpha}=\suml_{\NR{\alpha}=\n}\left(\sum^{+\infty}_{\k=-\infty}F_{\alpha}^{(\k)}\E^{\I\k\frac{2\pi}{\T}t}\right)\X^\alpha\nonumber\\
&=&\sum^{+\infty}_{\k=-\infty}\left(\suml_{\NR{\alpha}=\n}F_{\alpha}^{(\k)}\X^\alpha\right)\E^{\I\k\frac{2\pi}{\T}t}:=\sum^{+\infty}_{\k=-\infty} F^{(\k)}(\X)\E^{\I\k\frac{2\pi}{\T}t}.\nonumber
\end{eqnarray}
In this first step, we prove that if there exists $\Phi$ in $\Hom^{\n}(\Ezero,\Sob^{\l+1}(\R/\T\Z,\Eun))$ such that $\AAL\Phi=\F$, then necessarily $\Phi$ is unique. Indeed, for any function $\Phi$ in $\Hom^{\n}(\Ezero,\Sob^{\l+1}(\R/\T\Z,\Eun))$, we have
\begin{eqnarray}
\AAL\Phi=\F &\Longleftrightarrow & \NHP{\AAL\Phi-\F}{\n,\Sob^0}^2=0,\nonumber\\
										&\Longleftrightarrow & \sum^{+\infty}_{\k=-\infty}\NP{(\mathcal{A}_L+\I\k\frac{2\pi}{\T})\Phi^{(\k)}(\X)-\F^{(\k)}(\X)}{\n}^2=0,\nonumber\\
										&\Longleftrightarrow & \forall\k\in\Z, \quad(\mathcal{A}_L+\I\k\frac{2\pi}{\T})\Phi^{(\k)}(\X)-\F^{(\k)}(\X)=0. \nonumber
\end{eqnarray}
We will then use Lemma 21 b) in Iooss and Lombardi (2010), which will enable us to state that all the linear operators $(\mathcal{A}_L+\I\k\frac{2\pi}{\T})$, for $\k$ in $\Z$, are invertible in $\Hom^{\n}(\Ezero,\Eun)$. Indeed, for any polynomial $P$,
\begin{eqnarray}
(\mathcal{A}_L+\I\k\frac{2\pi}{\T})P(\X)&=&D_\X P(\X)L_0\X-L_1P(\X)+\I\k\frac{2\pi}{\T}P(\X)\nonumber\\
					&=&D_\X P(\X)L_0\X- (L_1-\I\k\frac{2\pi}{\T})P(\X)\nonumber\\
					&=&\mathcal{A}_{L^{(k)}}P(\X)\nonumber,
\end{eqnarray}
where 
\begin{equation}\label{defLk}
L^{(\k)}:=\left(
\begin{array}{cc}
L_0 & 0 \\
0 & L_1-\I\k\frac{2\pi}{\T}\Id
\end{array}
\right).
\end{equation}
Moreover, (\ref{diophantian}) in hypothesis \ref{hyp2} ensures that, for all $\k$, $L^{(\k)}$ satisfies the hypothesis of Lemma 21 in Iooss and Lombardi (2010). Hence, for all $\k$, $\mathcal{A}_{L^{(k)}}$ is invertible, and thus
\begin{eqnarray}
\AAL\Phi=\F &\Longleftrightarrow & \forall\k\in\Z, \quad\mathcal{A}_{L^{(\k)}}\Phi^{(\k)}(\X)-\F^{(\k)}(\X)=0, \nonumber\\
                        &\Longleftrightarrow & \forall\k\in\Z, \quad\Phi^{(\k)}(\X)=\mathcal{A}_{L^{(\k)}}^{-1}\F^{(\k)}(\X).\nonumber
\end{eqnarray}
Then, necessarily, if $\Phi$ exists, then
\begin{equation}\label{defphi}
\Phi(\X,t)=\sum^{+\infty}_{\k=-\infty} (\mathcal{A}_L+\I\k\frac{2\pi}{\T})^{-1}F^{(\k)}(\X)\E^{\I\k\frac{2\pi}{\T}t}.
\end{equation}
\vspace{2ex}

\textbf{Step 2: Now, we prove that $\AAn$ is also surjective.}

Our aim is now to prove that the $\Phi$ defined by (\ref{defphi}) is well-defined and belongs to $\Hom^{\n}(\Ezero,\Sob^{\l+1}(\R/\T\Z,\Eun))$.

\vspace{1ex}

\textbf{Step 2.1: $\Phi$ is well-defined.}

To show that $\Phi$ is well-defined, we prove that $\Phi$ is in $\Hom^{\n}(\Ezero,\Sob^{0}(\R/\T\Z,\Eun))$, proving that
$$\sum^{+\infty}_{\k=-\infty} \NP{(\mathcal{A}_L+\I\k\frac{2\pi}{\T})^{-1}F^{(\k)}(\X)}{\n}^2<+\infty.$$
We use, here again, inequality (\ref{diophantian}) and Lemma 21 in Iooss and Lombardi (2010), with $L^{(k)}$ defined in (\ref{defLk}), to check:
$$\NP{(\mathcal{A}_L+\I\k\frac{2\pi}{\T})^{-1}F^{(\k)}(\X)}{\n}\hspace{-1ex}=\NP{\mathcal{A}_{L^{(\k)}}^{-1}F^{(\k)}(\X)}{\n}\hspace{-1ex}\leq\frac{\nu}{\ga^{\nu}}(\n+\NR{\k})^{\ta\nu}\NP{F^{(\k)}(\X)}{\n}.$$
Then, since by hypothesis \ref{hyp2} $\ta\nu\leq\l$ holds, we obtain
\begin{eqnarray}
\sum^{+\infty}_{\k=-\infty} \NP{(\mathcal{A}_L+\I\k\frac{2\pi}{\T})^{-1}F^{(\k)}(\X)}{\n}^2\hspace{-1ex}& \leq &  \sum^{+\infty}_{\k=-\infty}\frac{\nu^2}{\ga^{2\nu}}(\n+\NR{\k})^{2\ta\nu}\NP{F^{(\k)}(\X)}{\n}^2 \nonumber \\
&\leq& \frac{\nu^2\n^{2\ta\nu}}{\ga^{2\nu}}\sum^{+\infty}_{\k=-\infty} (1+\k^2)^{\ta\nu}\NP{F^{(\k)}(\X)}{\n}^2\nonumber\\
&\leq& \frac{\nu^2\n^{2\ta\nu}}{\ga^{2\nu}}\sum^{+\infty}_{\k=-\infty} (1+\k^2)^{\l}\NP{F^{(\k)}(\X)}{\n}^2\nonumber\\
\NHP{\Phi}{n,\Sob^0}^2 &\leq& \frac{\nu^2\n^{2\ta\nu}}{\ga^{2\nu}} \NHP{F}{\n}^2.\nonumber
\end{eqnarray}
Thus $\Phi$ is well-defined and
\begin{equation}\label{init}
\NHP{\Phi}{n,\Sob^0} \leq \frac{\nu\n^{\ta\nu}}{\ga^{\nu}} \NHP{F}{\n}.
\end{equation}

\vspace{1ex}

\textbf{Step 2.2: $\Phi$ belongs to $\Hom^{\n}(\Ezero,\Sob^{\l+1}(\R/\T\Z,\Eun))$.}

We show, by induction, that if $\j\leq\l+1$, then
\begin{equation}\label{hered}
\Phi\in\Hom^{\n}(\Ezero,\Sob^{\j}(\R/\T\Z,\Eun))\text{ and }\NHP{\Phi}{\n,\Sob^\j} \leq \C_\j n^{\j+\ta\nu}\NHP{F}{\n}.
\end{equation}
First, observe that (\ref{init}) ensures that (\ref{hered}) holds for $\j=0$. Then, we assume that (\ref{hered}) holds for one $\j\leq\l$. We know that, for all $\k$ in $\Z$,
$$(\mathcal{A}_L+\I\k\frac{2\pi}{\T})\Phi^{(\k)}(\X)=F^{(\k)}(\X).$$
Hence,
$$\I\k\frac{2\pi}{\T}\Phi^{(\k)}(\X)=F^{(\k)}(\X)-\mathcal{A}_L\Phi^{(\k)}(\X).$$
So
$$(1+\k^2)\NP{\Phi^{(\k)}}{\n}^2=\NP{\Phi^{(\k)}}{\n}^2+\frac{\T^2}{4\pi^2}\NP{F^{(\k)}-\mathcal{A}_L\Phi^{(\k)}}{\n}^2,$$
and
\begin{equation}\label{ineq}
(1+\k^2)\NP{\Phi^{(\k)}}{\n}^2\leq\NP{\Phi^{(\k)}}{\n}^2+\frac{\T^2}{2\pi^2}\left(\NP{F^{(\k)}}{\n}^2+\NP{\mathcal{A}_L\Phi^{(\k)}}{\n}^2\right). 
\end{equation}
The proof of lemma 21 in Iooss and Lombardi (2010) ensures that the eigenvalues of $\mathcal{A}_L$ on 

$\Hom^{\n}(\Ezero,\Eun))$ are the
$$|\scal{\a}{\lambda^{(0)}}-\lambda_\j^{(1)}|,\quad\text{for }|\a|=\n, 1\leq\j\leq\m_1.$$
Then
\begin{eqnarray}
\NP{\mathcal{A}_L\Phi^{(\k)}(\X)}{\n}\hspace{-2ex}&\leq&\hspace{-2ex}\max_{|\a|=\n, 1\leq\j\leq\m_1}|<\a,\lambda^{(0)}>-\lambda_\j^{(1)}|\NP{\Phi^{(\k)}(\X)}{\n}\nonumber\\
&\leq&\hspace{-2ex} (\n+1) \max_{\dindice{1\leq\i\leq\m_0}{1\leq\j\leq\m_1}}\{|\lambda_\i^{(0)}|,|\lambda_\j^{(1)}|\}\NP{\Phi^{(\k)}(\X)}{\n}\nonumber\\ 
&=&(\n+1)\Lambda\NP{\Phi^{(\k)}(\X)}{\n}.\label{AL}
\end{eqnarray}
with the $\Lambda$ defined in the stating of the lemma. Then, combining (\ref{ineq}) and (\ref{AL}), we obtain that, for all $\k$ in $\Z$,
$$(1+\k^2)\NP{\Phi^{(\k)}}{\n}^2\leq(1+\frac{\T^2}{2\pi^2}\Lambda^2(\n+1)^2)\NP{\Phi^{(\k)}}{\n}^2+\frac{\T^2}{2\pi^2}\NP{F^{(\k)}}{\n}^2.$$
Hence, since we assume that (\ref{hered}) holds for $\j$,
\begin{eqnarray}
\NHP{\Phi}{n,\Sob^{\j+1}}^2\hspace{-2ex}&=&\hspace{-1ex}\sum^{+\infty}_{\k=-\infty}(1+\k^2)^{\j+1}\NP{\Phi^{(\k)}}{\n}^2\nonumber\\
&\leq&\hspace{-2ex}\sum^{+\infty}_{\k=-\infty}\hspace{-0,5ex}(1+\k^2)^{\j}\left(\hspace{-0,5ex}\left(1\hspace{-0,5ex}+\hspace{-0,5ex}\frac{1}{2}\left(\frac{\T\Lambda(\n+1)}{\pi}\right)^2\right)\NP{\Phi^{(\k)}}{\n}^2\hspace{-2ex}+\frac{\T^2}{2\pi^2}\NP{F^{(\k)}}{\n}^2\right)\nonumber\\
&\leq&\left(1+\frac{1}{2}\left(\frac{\T\Lambda(\n+1)}{\pi}\right)^2\right)\NHP{\Phi}{\n,\Sob^\j}^2+\frac{\T^2}{2\pi^2}\NHP{F}{\n,\Sob^j}^2\nonumber\\
&\leq&\left(\left(1+\frac{1}{2}\left(\frac{\T\Lambda(\n+1)}{\pi}\right)^2\right)\C_\j^2\n^{2\j+2\ta\nu} + \frac{\T^2}{2\pi^2}\right)\NHP{F}{\n,\Sob^\l}^2\nonumber\\
&\leq& \C_{\j+1}^2\n^{2(\j+1+\ta\nu)}\NHP{F}{\n}^2,\nonumber
\end{eqnarray}
which means that
\begin{equation}
\NHP{\Phi}{n,\Sob^{\j+1}} \leq \C_{\j+1} \n^{\j+1+\ta\nu}\NHP{F}{\n}.\nonumber
\end{equation}
%
%
%
\cqfd

\subsection{A few properties of norms}\label{propnorm}

We will need, in what follows, a few properties of multiplicativity for the new norms introduced in section \ref{defnorm}. They will be usefull to compute upper bounds for the remainder, whose expression is given with symmetric $q$-linear applications.
We gather these properties and their proofs in this section.

\begin{lem}\label{propnorm1}
If $\Phi$ belongs to $\Hom^{\n}(\Ezero,\Sob^{\l}(\R/\T\Z,\R^\m))$, then
$$\NH{\Phi(\uzero,.)}{\l}\leq\NHP{\Phi}{\n}\NR{\uzero}^\n.$$
\end{lem}

\textbf{Proof.}
If $\Phi$ is in $\Hom^{\n}(\Ezero,\Sob^{\l}(\R/\T\Z,\R^\m))$, it reads
$$\Phi(\X,t)=\sum_{|\alpha|=\n}\phi_\alpha(t)\X^{\alpha},$$
with $\phi_\alpha$ in $\Sob^\l(\R/\T\Z,\R^\m)$. Then, with $\uzero=x_1 e_1+\cdots+x_{\m_0}e_{\m_0}$ (notations of section \ref{defnorm}),
\begin{eqnarray}
\NH{\Phi(\uzero,.)}{\l}\hspace{-1ex}&=&\hspace{-1ex}\NH{\sum_{|\alpha|=\n}\phi_\alpha(t)\uzero^{\alpha}}{\l}\leq\sum_{|\alpha|=\n}\NH{\phi_\alpha}{\l}|x_1|^{\alpha_1}\cdots |x_{\m_0}|^{\alpha_{\m_0}}\nonumber\\
&=&\hspace{-1ex}(\sum_{|\alpha|=\n}\NH{\phi_\alpha}{\l}\X^\alpha)\Bigr |_{\X=(|x_1|,\cdots,|x_{\m_0}|)}\leq\Norme{\sum_{|\alpha|=\n}\NH{\phi_\alpha}{\l}\X^\alpha}{0,\n}\hspace{-2ex}\NR{\uzero}^\n,\nonumber
\end{eqnarray}
where $\Norme{.}{0,\n}$ is the norm defined on $\Hom^{\n}(\R^\m,\R^\m)$ by
$$\Norme{P(\X)}{0,\n}:=\supl_{u\in\R^\m}\frac{\NR{P(u)}}{\NR{u}^\m}.$$
Lemma 2.10 of Iooss and Lombardi (2005) ensures that, for all $P$ in $\Hom^{\n}(\R^\m,\R^\m)$, we have
$$\Norme{P(\X)}{0,\n}\leq\Norme{P(\X)}{2,\n}.$$
So,
\begin{equation}
\NH{\Phi(\uzero,.)}{\l}\leq\Norme{\sum_{|\alpha|=\n}\NH{\phi_\alpha}{\l}\X^\alpha}{2,\n}\NR{\uzero}^\n=\NHP{\Phi}{\n}\NR{\uzero}^\n.\nonumber
\end{equation}\cqfd

\vspace{1ex}

We then admit the following
\begin{lem}\label{lem2.2}
There exists a constant $\CSob$ such that for all functions $f$, $g$ in $\Sob^{\l}(\R/\T\Z,\R^\m)$,
$$\NH{fg}{\l}\leq\CSob\NH{f}{\l}\NH{g}{\l}.$$
\end{lem}

And we use it to prove the following

\begin{lem}\label{propnorm3}
Let $R_{\q} \left(t\right)\left[\X_1,\cdots,\X_\q\right]$ be a symmetric $\q$-linear application 

from $(\R^\m)^\q$ to $\Ezero$ or $\Eun$. Then, we have
$$(A_1) \Longrightarrow (A_2) \Longrightarrow (A_3) \Longrightarrow (B) \Longrightarrow (A_1'),$$
where
\begin{description}
\item[$(A_1)$] for all $x_1,\cdots,x_\q$ in $\R^\m$, $\NH{R_{\q}\left(t\right)\left[x_1,\cdots,x_\q\right]}{\l}\leq\c\left(\Frac{1}{\rh}\right)^\q\NR{x_1}\cdots\NR{x_\q};$
\item[$(A_2)$] denoting 
$$R_\q(t)\left[\X_1,\cdots,\X_\q\right]:=\sum_{1\leq i_\j\leq\m}a_{(i_1,\cdots,i_\q)}(t)\X_{1,i_1}\cdots\X_{\q,i_\q},$$
then, for all $x_1,\cdots,x_\q$ in $\R^\m$
$$\NR{\sum_{1\leq i_\j\leq\m}\NH{a_{(i_1,\cdots,i_\q)}}{\l}x_{1,i_1}\cdots x_{\q,i_\q}}\leq \c \left(\frac{\sqrt{\m}}{\rh}\right)^\q\NR{x_1}\cdots\NR{x_\q};$$
\item[$(A_3)$] for all $f_1,\cdots,f_\q$ in $\Sob^\l(\R/\T\Z,\R^\m)$, $$\NH{R_{\q}\left(t\right)\left[f_1,\cdots,f_\q\right]}{\l}\leq\c\left(\frac{\CSob\sqrt{\m}}{\rh}\right)^\q\NH{f_1}{\l}\cdots\NH{f_\q}{\l};$$
\item[$(B)$] if $\k_1+\cdots+\k_\q=\n$, then for all $\Phi_{\k_1},\cdots,\Phi_{\k_\q}$ in $\Hom^{\k_i}(\Ezero,\Sob^\l(\R/\T\Z,\R^\m))$
$$\NHP{R_\q(t)\left[\Phi_{\k_1},\cdots,\Phi_{\k_\q}\right]}{\n}\leq\c\left(\frac{\CSob\sqrt{\m}}{\rh}\right)^\q\NHP{\Phi_{\k_1}}{\k_1}\cdots\NHP{\Phi_{\k_\q}}{\k_\q};$$
\item[$(A_1')$] \hspace{-1ex}for all $x_1,\cdots,x_\q$ in $\R^\m$\hspace{-0,5ex}, $\NH{R_{\q}\left(t\right)\left[x_1,\cdots,x_\q\right]}{\l}\hspace{-1ex}\leq\c\left(\Frac{\CSob\sqrt{\m}}{\rh}\right)^\q\NR{x_1}\cdots\NR{x_\q}.$
\end{description}
\end{lem}
\begin{rem}
$(A_1)$ and $(A_1')$ are nearly the same property: the only difference between them is that in $(A_1')$ the constant in the upper bound is bigger than in $(A_1)$. This is due to the fact that in each implication we lose precision in the upper bounds.
\end{rem}
\begin{rem}
For the following, the most usefull results in this lemma are $(A_1) \Rightarrow (A_3)$ and $(A_1) \Rightarrow (B)$.
\end{rem}

\textbf{Proof of $(A_1) \Rightarrow (A_2)$.}

Assume that $(A_1)$ holds. Let $(e_1,\cdots,e_\m)$ be the basis of $\R^\m$ introduced in section \ref{defnorm}. Take $x_1=e_{i_1},\cdots,x_\q=e_{i_\q}$, then
$$\NH{R_{\q}\left(t\right)\left[e_{i_1},\cdots,e_{i_\q}\right]}{\l}=\NH{a_{(i_1,\cdots,i_\q)}(t)}{\l}\leq\c\left(\frac{1}{\rh}\right)^\q\NR{e_1}\cdots\NR{e_{i_\q}}=\frac{\c}{\rh^\q}.$$
Hence, by Cauchy-Schwarz inequality,
\begin{eqnarray}
&&\hspace{-10ex}\NR{\sum_{1\leq i_1,\cdots,i_\q\leq\m}\NH{a_{(i_1,\cdots,i_\q)}}{\l}x_{1,i_1}\cdots x_{\q,i_\q}}\leq\frac{\c}{\rh^\q}\NR{\sum_{1\leq i_1,\cdots,i_\q\leq\m}|x_{1,i_1}\cdots x_{\q,i_\q}|}\nonumber\\
&=&\frac{\c}{\rh^\q}\left(\sum_{i_1=1}^{\m}|x_{1,i_1}|\right)\times\cdots\times\left(\sum_{i_\q=1}^{\m}|x_{\q,i_\q}|\right)\nonumber\\
&\leq&\frac{\c}{\rh^\q}\left(\sqrt{\m}\sqrt{\sum_{i_1=1}^{\m}|x_{1,i_1}|^2}\right)\times\cdots\times\left(\sqrt{\m}\sqrt{\sum_{i_\q=1}^{\m}|x_{\q,i_\q}|^2}\right)\nonumber\\
&=&\c\left(\frac{\sqrt{\m}}{\rh}\right)^\q\NR{x_1}\cdots\NR{x_\q}.\nonumber
\end{eqnarray}

\textbf{Proof of $(A_2) \Rightarrow (A_3)$.}

Assume that $(A_2)$ holds. Let $f_1,\cdots,f_\q$ be functions of $\Sob^\l(\R/\T\Z,\R^\m)$. Then, using lemma \ref{lem2.2},
\begin{eqnarray}
\NH{R_{\q}\left(t\right)\left[f_1,\cdots,f_\q\right]}{\l}&=&\NH{\sum_{1\leq i_1,\cdots,i_\q\leq\m}a_{(i_1,\cdots,i_\q)}(t)f_{1,i_1}(t)\cdots f_{\q,i_\q}(t)}{\l}\nonumber\\
		&\leq&\sum_{1\leq i_1,\cdots,i_\q\leq\m}\CSob^\q\NH{a_{(i_1,\cdots,i_\q)}}{\l}\NH{f_{1,i_1}}{\l}\cdots \NH{f_{\q,i_\q}(t)}{\l}\nonumber\\
		&\leq&\CSob^\q \c\left(\frac{\sqrt{\m}}{\rh}\right)^\q\NH{f_1}{\l}\cdots\NH{f_\q}{\l}\nonumber\\
		&=&\c\left(\frac{\CSob\sqrt{\m}}{\rh}\right)^\q\NH{f_1}{\l}\cdots\NH{f_\q}{\l}.\nonumber
\end{eqnarray}

\textbf{Proof of $(A_3)\Rightarrow (B)$.}

We follow the same strategy of proof as that of lemma A8 of Iooss and Lombardi (2005). Assume that $(A_3)$ is satisfied. Take $\k_1,\cdots,\k_\q$ such that 
$$\k_1+\cdots+\k_\q=\n,$$
and take $\Phi_{\k_i}$ in $\Hom^{\k_i}(\Ezero,\Sob^\l(\R/\T\Z,\R^\m))$, for $1\leq i\leq\q$. Denote
$$\Phi_{\k_i}(\X,t):=\sum_{|\alpha|=\k_i}\phi_{\alpha,\k_i}(t)\X^\alpha.$$
Then,
\begin{eqnarray}
R_\q(t)\left[\Phi_{\k_1},\cdots,\Phi_{\k_\q}\right]&=&\hspace{-1ex}\sum_{\dindice{\alpha^{(\i)}\in\N^{\m_0}}{|\alpha^{(\i)}|=\k_\i}}\X^{\alpha^{(1)}+\cdots +\alpha^{(\q)}}R_\q(t)\left[\phi_{\alpha^{(1)},\k_1},\cdots ,\phi_{\alpha^{(\q)},\k_\q}\right]\nonumber\\
&=&\sum_{|\be|=\n}\X^\be\hspace{-3ex}\sum_{\dindice{|\alpha^{(\i)}|=\k_\i}{\alpha^{(1)}+\cdots +\alpha^{(\q)}=\be}}\hspace{-3ex}R_\q(t)\left[\phi_{\alpha^{(1)},\k_1},\cdots ,\phi_{\alpha^{(\q)},\k_\q}\right].\nonumber
\end{eqnarray}
Hence, using Cauchy-Schwarz inequality again,
\begin{eqnarray}
&&\hspace{-4ex}\NHP{R_\q(t)\left[\Phi_{\k_1},\cdots,\Phi_{\k_\q}\right]}{\n}^2=\sum_{|\be|=\n}\frac{\be!}{\n!}\NH{\hspace{-4ex}\sum_{\dindice{|\alpha^{(\i)}|=\k_\i}{\hspace{4ex}\alpha^{(1)}+\cdots +\alpha^{(\q)}=\be}}\hspace{-5ex}R_\q(t)\left[\phi_{\alpha^{(1)},\k_1},\cdots ,\phi_{\alpha^{(\q)},\k_\q}\right]}{\l}^2\nonumber
\end{eqnarray}
\vspace{-5ex}
\begin{eqnarray}
&\leq&\sum_{|\be|=\n}\frac{\be!}{\n!}\left(\c\left(\frac{\CSob\sqrt{\m}}{\rh}\right)^\q\sum_{\dindice{\alpha^{(1)}+\cdots +\alpha^{(\q)}=\be}{\alpha^{(\i)}\in\N^{\m_0},|\alpha^{(\i)}|=\k_\i}}\NH{\phi_{\alpha^{(1)},\k_1}}{\l}\cdots \NH{\phi_{\alpha^{(\q)},\k_\q}}{\l} \right)^2\nonumber\\
&\leq&\sum_{|\be|=\n}\frac{\be!}{\n!}\c^2\left(\frac{\CSob\sqrt{\m}}{\rh}\right)^{2\q}
\left(\sum_{\dindice{\alpha^{(1)}+\cdots +\alpha^{(\q)}=\be}{\alpha^{(\i)}\in\N^{\m_0},|\alpha^{(\i)}|=\k_\i}}\frac{1}{\alpha^{(1)}!}\cdots \frac{1}{\alpha^{(\q)}!}\right)\nonumber\\
&&\hspace{5ex}\times\left(\sum_{\dindice{\alpha^{(1)}+\cdots +\alpha^{(\q)}=\be}{\alpha^{(\i)}\in\N^{\m_0},|\alpha^{(\i)}|=\k_\i}}\left(\alpha^{(1)}!\NH{\phi_{\alpha^{(1)},\k_1}}{\l}^2\right)\cdots \left(\alpha^{(\q)}!\NH{\phi_{\alpha^{(\q)},\k_\q}}{\l}^2\right)\right).\nonumber
\end{eqnarray}
And lemma A6 of Iooss and Lombardi (2005) ensures that
$$\sum_{\dindice{\alpha^{(1)}+\cdots +\alpha^{(\q)}=\be}{\alpha^{(\i)}\in\N^{\m_0},|\alpha^{(\i)}|=\k_\i}}\frac{1}{\alpha^{(1)}!}\cdots \frac{1}{\alpha^{(\q)}!}=\frac{\n!}{\be!}\frac{1}{\k_1!\cdots \k_\q!}.$$
Then 
\begin{eqnarray}
&&\hspace{-6ex}\frac{1}{\c^2}\left(\frac{\rh}{\CSob\sqrt{\m}}\right)^{2\q}\NHP{R_\q(t)\left[\Phi_{\k_1},\cdots,\Phi_{\k_\q}\right]}{\n}^2\nonumber\\
&\leq&\sum_{|\be|=\n}\frac{1}{\k_1!\cdots \k_\q!}
\sum_{\dindice{\alpha^{(1)}+\cdots +\alpha^{(\q)}=\be}{|\alpha^{(\i)}|=\k_\i}}\left(\alpha^{(1)}!\NH{\phi_{\alpha^{(1)},\k_1}}{\l}^2\right)\cdots \hspace{-1ex}\left(\alpha^{(\q)}!\NH{\phi_{\alpha^{(\q)},\k_\q}}{\l}^2\right)\nonumber\\
&=&\sum_{\dindice{\alpha^{(\i)}\in\N^{\m_0}}{|\alpha^{(\i)}|=\k_\i}}\left(\frac{\alpha^{(1)}!}{\k_1!}\NH{\phi_{\alpha^{(1)},\k_1}}{\l}^2\right)\cdots \left(\frac{\alpha^{(\q)}!}{\k_\q!}\NH{\phi_{\alpha^{(\q)},\k_\q}}{\l}^2\right)\nonumber\\
&=&\left(\sum_{|\alpha^{(1)}|=\k_1}\left(\frac{\alpha^{(1)}!}{\k_1!}\NH{\phi_{\alpha^{(1)},\k_1}}{\l}^2\right)\right)\cdots\left(\sum_{|\alpha^{(1)}|=\k_\q}\left(\frac{\alpha^{(\q)}!}{\k_\q!}\NH{\phi_{\alpha^{(\q)},\k_\q}}{\l}^2\right)\right)\nonumber\\
&=&\NHP{\Phi_{\k_1}}{\k_1}^2\cdots\NHP{\Phi_{\k_\q}}{\k_\q}^2.\nonumber
\end{eqnarray}
Finally, we have
$$\NHP{R_\q(t)\left[\Phi_{\k_1},\cdots,\Phi_{\k_\q}\right]}{\n}\leq\c\left(\frac{\CSob\sqrt{\m}}{\rh}\right)^\q\NHP{\Phi_{\k_1}}{\k_1}\cdots\NHP{\Phi_{\k_\q}}{\k_\q}.$$

\textbf{Proof of $(B) \Rightarrow (A_1').$}

For given $x_1,\cdots,x_\q$ in $\R^\m$, take
$$\k_1=\cdots=\k_\q=\n=0, \quad\Phi_{\k_1}(\X,t)=x_1, \cdots, \Phi_{\k_q}(\X,t)=x_\q.$$
Then, $(B)$ ensures that
\begin{eqnarray}
\NH{R_{\q}\left(t\right)\left[x_1,\cdots,x_\q\right]}{\l}\hspace{-2ex}&=&\NHP{R_{\q}\left(t\right)\left[x_1,\cdots,x_\q\right]}{0}\leq\c\left(\frac{\CSob\sqrt{\m}}{\rh}\right)^\q\NHP{x_1}{0}\cdots\NHP{x_\q}{0}\nonumber\\
		&=&\hspace{-1,5ex}\c\left(\frac{\CSob\sqrt{\m}}{\rh}\right)^\q\hspace{-0,5ex}\NH{x_1}{\l}\cdots\NH{x_\q}{\l}=\c\left(\frac{\CSob\sqrt{\m}}{\rh}\right)^\q\hspace{-0,5ex}\NR{x_1}\cdots\NR{x_\q}.\nonumber
\end{eqnarray}\cqfd

\begin{lem}\label{propnorm4}
If $\Phi_{\k}$ is in $\Hom^{\k}(\Ezero,\Sob^\l(\R/\T\Z,\R^\m))$ and $\NN_\p$ in 

$\Hom^{\p}(\Ezero,\Sob^\l(\R/\T\Z,\Ezero))$, then
$$\NHP{D_\X\Phi_\k(\X,t).\NN_\p(\X,t)}{\k-1+\p}\leq\CSob\k\sqrt{\m_0}\NHP{\Phi_\k}{\k}\NHP{\NN_\p}{\p}.$$
\end{lem}

\textbf{Proof.} Denote
\begin{equation}
\Phi_\k(\X,t)=\sum_{|\alpha|=\k}\phi_\alpha(t)\X^\alpha, \quad \NN_\p(\X,t)=\sum_{|\be|=\p}\Nn_\be(t)\X^\be,\nonumber
\end{equation}
where the $\phi_\alpha$ are in $\Sob^\l(\R/\T\Z,\R^\m)$ and $\Nn_\be$ in $\Sob^\l(\R/\T\Z,\Ezero)$. Moreover, we denote
$$\Nn_\be(t):=(\Nn_{\be,1}(t),\cdots,\Nn_{\be,\m_0}(t)).$$
Then, we have
\begin{eqnarray}
D_\X\Phi_\k(\X,t).\NN_\p(\X,t)&=&\sum_{\j=1}^{\m_0}\bigg(\sum_{\dindice{|\alpha|=\k}{|\be|=\p}}\phi_\alpha(t)\alpha_\j\X^{\alpha-\sig_\j}\Nn_{\be,\j}(t)\X^\be\bigg)\nonumber\\
&=&\sum_{|\gam|=\k-1+\p}\X^\gam\bigg(\sum_{\j=1}^{\m_0}\sum_{\dindice{|\alpha|=\k,|\be|=\p}{\alpha+\be-\sig_\j=\gam}}\alpha_\j\phi_\alpha(t)\Nn_{\be,\j}(t)\bigg),\nonumber
\end{eqnarray}
where $\sig_\j$ stands for $(0,\cdots,0,\overbrace{1}^{\j^{\text{th}}},0,\cdots,0)$ in $\N^{\m_0}$. Thus, using Cauchy-Schwarz inequality, we get
\begin{eqnarray}
&&\hspace{-5ex}\NHP{D_\X\Phi_\k(\X,t).\NN_\p(\X,t)}{\k-1+\p}^2\nonumber\\
&=&\sum_{|\gam|=\k-1+\p}\frac{\gam!}{(\k-1+\p)!}\NH{\sum_{\j=1}^{\m_0}\sum_{\dindice{|\alpha|=\k,|\be|=\p}{\alpha+\be-\sig_\j=\gam}}\alpha_\j\phi_\alpha(t)\Nn_{\be,\j}(t)}{\l}^2\nonumber\\
&\leq&\sum_{|\gam|=\k-1+\p}\frac{\gam!}{(\k-1+\p)!}\bigg(\sum_{\j=1}^{\m_0}\sum_{\dindice{|\alpha|=\k,|\be|=\p}{\alpha+\be-\sig_\j=\gam}}|\alpha_\j|\CSob\NH{\phi_\alpha}{\l}\NH{\Nn_{\be,\j}}{\l} \bigg)^2\nonumber\\
&\leq&\CSob^2\sum_{|\gam|=\k-1+\p}\frac{\gam!}{(\k-1+\p)!}\bigg(\sum_{\j=1}^{\m_0}\sum_{\dindice{|\alpha|=\k,|\be|=\p}{\alpha+\be-\sig_\j=\gam}}\frac{\alpha_\j^2}{\alpha!\be!}\bigg)\nonumber\\
&&\hspace{20ex}\times\bigg(\sum_{\j=1}^{\m_0}\sum_{\dindice{|\alpha|=\k,|\be|=\p}{\alpha+\be-\sig_\j=\gam}}\alpha!\be!\NH{\phi_\alpha}{\l}^2\NH{\Nn_{\be,\j}}{\l}^2\bigg).
\nonumber
\end{eqnarray}
Since lemma A7 of Iooss and Lombardi (2005) ensures that
$$\sum_{\j=1}^{\m_0}\sum_{\dindice{|\alpha|=\k,|\be|=\p}{\alpha+\be-\sig_\j=\gam}}\frac{\alpha_\j^2}{\alpha!\be!}=\frac{1}{\k!\p !}(\k^2+(\m_0-1)\k)\frac{(\k-1+\p)!}{\gam!},$$
we get
\begin{eqnarray}
&&\hspace{-7ex}\NHP{D_\X\Phi_\k(\X,t).\NN_\p(\X,t)}{\k-1+\p}^2\nonumber\\
&\leq&\CSob^2(\k^2+(\m_0-1)\k)\sum_{|\gam|=\k-1+\p}\sum_{\j=1}^{\m_0}\sum_{\dindice{|\alpha|=\k,|\be|=\p}{\alpha+\be-\sig_\j=\gam}}\frac{\alpha!\be!}{\k!\p!}\NH{\phi_\alpha}{\l}^2\NH{\Nn_{\be,\j}}{\l}^2\nonumber\\
&=&\CSob^2(\k^2+(\m_0-1)\k)\sum_{|\alpha|=\k,|\be|=\p}\left(\frac{\alpha!}{\k!}\NH{\phi_\alpha}{\l}^2\left(\sum_{\j=1}^{\m_0}\frac{\be!}{\p!}\NH{\Nn_{\be,\j}}{\l}^2\right)\right)\nonumber\\
&=&\CSob^2(\k^2+(\m_0-1)\k)\left(\sum_{|\alpha|=\k}\frac{\alpha!}{\k!}\NH{\phi_\alpha}{\l}^2\right)\left(\sum_{|\be|=\p}\frac{\be!}{\p!}\NH{\Nn_\be}{\l}^2\right)\nonumber\\
&=&\CSob^2(\k^2+(\m_0-1)\k)\NHP{\Phi_\k}{\k}^2\NHP{\NN_\p}{\p}^2\nonumber\\
&\leq&\CSob^2\k^2\m_0\NHP{\Phi_\k}{\k}^2\NHP{\NN_\p}{\p}^2.\nonumber
\end{eqnarray}
Then, we finally have
$$\NHP{D_\X\Phi_\k(\X,t).\NN_\p(\X,t)}{\k-1+\p}\leq\CSob\k\sqrt{\m_0}\NHP{\Phi_\k}{\k}\NHP{\NN_\p}{\p}.$$\cqfd

\subsection{Step B: Choice of $\p=\p_{opt}$, upper bound for the remainder}\label{upperbound}
As described in part \ref{sketch}, for a fixed integer $\p$, with the $\Phi$ constructed in section \ref{construct}, we set
\begin{eqnarray}
\Reste(\uzero,t)&:=&(\Id-\Pi_{\p})(\V_1(\uzero+\Phi(\uzero,t),t) \label{reste}\\
										&&\hspace{5ex}-D_{\uzero}\Phi(\uzero,t).\V_0(\uzero+\Phi(\uzero,t),t)),\nonumber
\end{eqnarray}
and we compute an upper bound for $\NH{\Reste(\uzero,.)}{\l}$. To simplify notations, let us denote
\begin{equation}\label{defPPhi}
\Phi_1(\uzero,t):=\uzero,\quad\PPhi(\uzero,t):=\sum_{\n=1}^{\p}\Phi_\n(\uzero,t).
\end{equation}
Then, (\ref{reste}) reads
\begin{eqnarray}
\Reste(\uzero,t)&=&\sum_{\q=2}^{\p}\hspace{-2ex}\sum_{\dindice{1\leq\k_\j\leq\p}{\k_1+\cdots+\k_\q\geq\p+1}}\V_{1,\q}(t)[\Phi_{\k_1}(\uzero,t),\cdots,\Phi_{\k_\q}(\uzero,t)]\nonumber\\
&& - \sum_{2\leq\j,\q\leq\p}\hspace{-3ex}\sum_{\dindice{1\leq\k_\j\leq\p}{\k_1+\cdots+\k_\q\geq\p-\j+2}}\hspace{-5ex}D_{\uzero}\Phi_\j(\uzero,t).\V_{0,\q}(t)[\Phi_{\k_1}(\uzero,t),\cdots,\Phi_{\k_\q}(\uzero,t)]\nonumber\\
&& + \sum_{\q\geq\p+1}\V_{1,\q}(t)[\PPhi(\uzero,t)^{(\q)}] \nonumber\\
&&- \sum_{\j=2}^{\p}\sum_{\q\geq\p+1}D_{\uzero}\Phi_\j(\uzero,t).\V_{0,\q}(t)[\PPhi(\uzero,t)^{(\q)}].\nonumber
\end{eqnarray}
So, to evaluate $\NH{\Reste(\uzero,.)}{\l}$, the first step is to compute upper bounds for the $\NH{\Phi_\n(\uzero,.)}{\l}$ constructed in Step A. And in fact, using lemma \ref{propnorm1}, it will be sufficient to compute upper bounds for the $\NHP{\Phi_\n}{\n}$.
Denote them by
\begin{equation}
\pphi_\n:=\NHP{\Phi_\n}{\n}, \quad 1\leq\n\leq\p.\nonumber
\end{equation}

\subsubsection{Upper bounds for the $\pphi_\n$}
To compute upper bounds for the $\pphi_\n$, we have to get back to the construction of the $\Phi_\n$: they were constructed by induction with the equation (\ref{projn}). This equation reads, explicitly
\begin{eqnarray}
(\mathcal{A}_L+\partial_t)\Phi_\n&=&\sum_{q=2}^{\n}\sum_{\dindice{1\leq\k_\i}{\k_1+\cdots+\k_\q=\n}}\V_{1,\q}[\Phi_{\k_1},\cdots,\Phi_{\k_\q}]\nonumber\\
&&- \sum_{\j=2}^{\n-1}\sum_{\q=2}^{\n-\j+1}\hspace{-2ex}\sum_{\dindice{1\leq\k_\i}{\k_1+\cdots+\k_\q=\n-\j+1}}\hspace{-2ex}D_\X\Phi_\j.\V_{0,\q}[\Phi_{\k_1},\cdots,\Phi_{\k_\q}].\nonumber
\end{eqnarray}
With this last equation, using inequality (\ref{majinv}) of lemma \ref{leminv}, $(A_1)\Rightarrow(B)$ of lemma \ref{propnorm3} (since $\V$ is in $\mathcal{A}(\Nu,\Sob^\l(\R/\T\Z, \R^\m)$ and thus satisfies $(A_1)$), and using lemma \ref{propnorm4}, we obtain
\begin{eqnarray}
\pphi_\n&\leq&\C_\l\n^{\l+\ta\nu}\Bigg( \sum_{\q=2}^{\n}\sum_{\dindice{\k_1+\cdots+\k_\q=\n}{1\leq\k_\i}}\c\left(\frac{\CSob\sqrt{\m}}{\rh}\right)^\q\pphi_{\k_1}\cdots\pphi_{\k_\q}\nonumber\\
& & +\sum_{\j=2}^{\n}\sum_{\q=2}^{\n-\j+1}\sum_{\dindice{\k_1+\cdots+\k_\q=\n-\j+1}{1\leq\k_\i}}\CSob\j\sqrt{\m_0}\pphi_\j\c\left(\frac{\CSob\sqrt{\m}}{\rh}\right)^\q\pphi_{\k_1}\cdots\pphi_{\k_\q} \Bigg).\nonumber
\end{eqnarray}
Thus,
\begin{eqnarray}
\hspace{-7ex}\pphi_\n&\leq&\max(1,\CSob)\C_\l\n^{\l+\ta\nu}\c\Bigg( \sum_{\q=2}^{\n}\sum_{\dindice{\k_1+\cdots+\k_\q=\n}{1\leq\k_\i}}\left(\frac{\CSob\sqrt{\m}}{\rh}\right)^\q\pphi_{\k_1}\cdots\pphi_{\k_\q}\nonumber\\ &&+\sum_{\j=2}^{\n}\sum_{\q=2}^{\n-\j+1}\hspace{-3ex}\sum_{\dindice{1\leq\k_\i}{\k_1+\cdots+\k_\q=\n-\j+1}}\hspace{-3ex}\j\sqrt{\m_0}\left(\frac{\CSob\sqrt{\m}}{\rh}\right)^\q\pphi_\j\pphi_{\k_1}\cdots\pphi_{\k_\q}\Bigg) .\label{ineg}
\end{eqnarray}

\begin{lem}\label{lemK}
For all $\n$ with $1\leq\n\leq\p$,
$$\pphi_\n\leq\sqrt{\m_0}K^{\n-1}(\n!)^{\l+1+\ta\nu},$$
where
\begin{equation}
K:=\max(\frac{9\CSob\sqrt{\m_0}\sqrt{\m}}{\rh},\frac{8\C_\l\c(\CSob\sqrt{\m_0})^3\sqrt{\m}}{\rh^2}).\nonumber
\end{equation}
\end{lem}

\textbf{Proof.} The same strategy of proof as that of lemma 14 in Iooss and Lombardi (2010) works: indeed, with our norm we also have $$\pphi_1=\sqrt{\m_0};$$ and our inequality (\ref{ineg}) is the same as inequality (19) in the proof of lemma 14, if one equates their $\rho$ with, in our notations, $\frac{\rh}{\CSob\sqrt{\m}}$, and equates $a:=\CSob\C_\l$ and $\tau':=\l+\ta\nu$.
Then, all the following computations of Iooss and Lombardi (2010) work similarly.

\cqfd

\subsubsection{Choice of $\p_{opt}$, upper bounds for $\PPhi$}

Since now, we fix $\delta>0$ and choose a value $\p_{opt}(\delta)$ for $\p$, and in this subsection we find an upper bound for $\PPhi$ (defined in (\ref{defPPhi})) with this value of $\p$.
\begin{lem}\label{popt}
Fix $\delta>0$. Denote
\begin{equation}
\b:=\TFrac{1}{1+\l+\ta\nu},\quad\p_{opt}:=\left[\TFrac{1}{(2\delta K)^{\b}}\right]\nonumber
\end{equation}
Let us choose $\p=\p_{opt}$. Then, we have, for $\NR{\uzero}\leq\delta$,
\begin{equation}
\NH{\PPhi(\uzero,.)}{\l}\leq2\NR{\uzero}\sqrt{\m_0}\leq2\delta\sqrt{\m_0}.\nonumber
\end{equation}
\end{lem}

\textbf{Proof.} Take $\p=\p_{opt}$ and $\NR{\uzero}\leq\delta$. Then, using lemmas \ref{propnorm1} and \ref{lemK}:
\begin{eqnarray}
\NH{\PPhi(\uzero,.)}{\l}&=&\NH{\sum_{\n=1}^{\p}\Phi_\n(\uzero,.)}{\l}\leq\sum_{\n=1}^{\p}\NH{\Phi_\n(\uzero,.)}{\l}\leq\sum_{\n=1}^{\p}\NHP{\Phi_\n}{\n}\NR{\uzero}^\n\nonumber\\
&\leq&\NR{\uzero}\sum_{\n=1}^{\p}\pphi_\n\delta^{\n-1}\leq\NR{\uzero}\sqrt{\m_0}\sum_{\n=1}^{\p}K^{\n-1}\delta^{\n-1}(\n!)^\b\nonumber\\
					&\leq&\NR{\uzero}\sqrt{\m_0}\sum_{\n=1}^{\p}(K\delta)^{\n-1}(\p^{\n-1})^\b=\NR{\uzero}\sqrt{\m_0}\sum_{\n=1}^{\p}(K\delta\p^b)^{\n-1}\nonumber\\
					&\leq&\NR{\uzero}\sqrt{\m_0}\sum_{\n=1}^{\p}\left(\frac{1}{2}\right)^{\n-1}\leq2\NR{\uzero}\sqrt{\m_0}.\nonumber
\end{eqnarray}\cqfd

\subsubsection{Upper bound for the remainder}

\begin{lem}
There exists $\delta_0$ such that, if $\delta<\delta_0$, then, with the choice $\p=\p_{opt}$ of lemma \ref{popt}, the remainder $\Reste$ satisfies
\begin{equation}
\supl_{\NR{\uzero}\leq\delta}\NH{\Reste(\uzero,.)}{\l}\leq\M\E^{-\frac{\om}{\delta^\b}},\nonumber
\end{equation}
with 
$$\M:=\c\left(\frac{73}{72}+\CSob\sqrt{\m_0}(2\sqrt{\m}+\frac{\sqrt{\m_0}}{72})\right),\quad \om:=\frac{\ln(2)}{2(2K)^\b},$$
\begin{equation}\label{defdelta0}
\delta_0:=min\left(\frac{1}{2K(2\E)^\b},\frac{\rh}{4\CSob\sqrt{\m}\sqrt{\m_0}}\right).
\end{equation}
\end{lem}

\textbf{Proof.} Fix $\uzero$ such that
$$\NR{\uzero}\leq\delta.$$
We have
\begin{eqnarray}
\Reste(\uzero,t)&=&\sum_{\q=2}^{\p}\hspace{-2ex}\sum_{\dindice{1\leq\k_\j\leq\p}{\k_1+\cdots+\k_\q\geq\p+1}}\hspace{-2ex}\V_{1,\q}(t)[\Phi_{\k_1}(\uzero,t),\cdots,\Phi_{\k_\q}(\uzero,t)]\nonumber\\
&& + \sum_{\q\geq\p+1}\V_{1,\q}(t)[\PPhi(\uzero,t)^{(\q)}]\nonumber\\
&& - \sum_{2\leq\j,\q\leq\p}\hspace{-4ex}\sum_{\dindice{1\leq\k_\j\leq\p}{\k_1+\cdots+\k_\q\geq\p-\j+2}}\hspace{-4ex}D_{\uzero}\Phi_\j(\uzero,t).\V_{0,\q}(t)[\Phi_{\k_1}(\uzero,t),\cdots,\Phi_{\k_\q}(\uzero,t)]\nonumber\\
&& - \sum_{\j=2}^{\p}\sum_{\q\geq\p+1}D_{\uzero}\Phi_\j(\uzero,t).\V_{0,\q}(t)[\PPhi(\uzero,t)^{(\q)}].\nonumber
\end{eqnarray}
One can check, with the same computations as those of Iooss and Lombardi (2010) in their proof of their lemma 16, that combining our lemma \ref{propnorm1} with $(A_1)\Rightarrow(A_3)$ of lemma \ref{propnorm3} and lemma \ref{lemK}, we get $$\hspace{-9ex}\NH{\sum_{\q=2}^{\p}\sum_{\dindice{\k_1+\cdots+\k_\q\geq\p+1}{1\leq\k_\j\leq\p}}\V_{1,\q}(.)[\Phi_{\k_1}(\uzero,.),\cdots,\Phi_{\k_\q}(\uzero,.)]}{\l}\leq\frac{\c}{72}\E^{-\frac{\ln2}{(2\delta K)^\b}};$$
and that, if $\delta\leq\frac{1}{2K(2\E)^\b}$, then, by lemmas \ref{propnorm1}, \ref{propnorm3} (with $(A_1)\Rightarrow(B)$) and \ref{lemK}, we obtain
$$\NH{\sum_{2\leq\j,\q\leq\p}\hspace{-4ex}\sum_{\dindice{1\leq\k_\j\leq\p}{\k_1+\cdots\hspace{-1ex}+\k_\q\geq\p-\j+2}}\hspace{-6ex}D_{\uzero}\Phi_\j(\uzero,.).\V_{0,\q}(.)[\Phi_{\k_1}(\uzero,.),\cdots,\Phi_{\k_\q}(\uzero,.)]}{\l}\hspace{-2,5ex}\leq\frac{\c\CSob\m_0}{72}\E^{-\frac{\ln2}{2(2\delta K)^\b}},$$
$$\hspace{-7ex}\NH{\sum_{\j=2}^{\p}\sum_{\q\geq\p+1}D_{\uzero}\Phi_\j(\uzero,t).\V_{0,\q}(t)[\PPhi(\uzero,t)^{(\q)}]}{\l}\leq2\CSob\c\sqrt{\m}\sqrt{\m_0}\E^{-\frac{\ln2}{2(2\delta K)^\b}};$$
and that, if $\delta\leq\frac{\rh}{4\CSob\sqrt{\m}\sqrt{\m_0}}$, then, using $(A_1)\Rightarrow(A_3)$ of lemma \ref{propnorm3}, we get 
$$\hspace{-35ex}\NH{\sum_{\q\geq\p+1}\V_{1,\q}(t)[\PPhi(\uzero,t)^{(\q)}]}{\l}\leq\c \E^{-\frac{\ln2}{(2\delta K)^\b}}.$$\cqfd

\subsection{Step C: Upper bound for $\V^1$}\label{StepC}

To complete the proof of the theorem \ref{th1}, it remains to show that (\ref{majvun}) holds for the choices we have done for $\p$ and $\Phi$. As said in part \ref{sketch}, we set 
\begin{eqnarray}
\V^1(\uzero,\vun,t)&:=&-\Reste(\uzero,t)+\V_1(\uzero+\vun+\Phi(\uzero,t),t)\nonumber\\
										&&-D_{\uzero}\Phi(\uzero,t).\V_0(\uzero+\vun+\Phi(\uzero,t),t)\nonumber\\
										&&-(\mathcal{A}_L+\partial_t)\Phi(\uzero,t).\nonumber
\end{eqnarray}
Combining it with (\ref{step1}) and (\ref{step2}), we get
\begin{eqnarray}
\V^1(\uzero,\vun,t)\hspace{-1ex}
										&=&\hspace{-1ex}\V_1(\uzero+\vun+\Phi(\uzero,t),t)-\V_1(\uzero+\Phi(\uzero,t),t)\nonumber\\
										&&\hspace{-2ex}-D_{\uzero}\Phi(\uzero,t)\big(\V_0(\uzero+\vun+\Phi(\uzero,t),t)-\V_0(\uzero+\Phi(\uzero,t),t)\big)\nonumber
\end{eqnarray}
\vspace{-4ex}
$$=\left((\Id-P_0)\hspace{-1ex}-D_{\uzero}\Phi(\uzero,t).P_0 \right)\left(\V(\uzero+\vun+\Phi(\uzero,t),t)\hspace{-1ex}-\V(\uzero+\Phi(\uzero,t),t)\right)\nonumber$$
where $P_0$ stands for the projection of $\R^\m$ on $\Ezero$.
We begin by showing an upper bound for $D_{\uzero}\Phi$.
\begin{lem}\label{lem1}
If $0<\delta<\delta_0$, and $\p=\p_{opt}$, then, for all $\F$ in $\Sob^\l(\R/\T\Z,\Ezero)$, and all $\NR{\uzero}\leq\delta$, we have
$$\NH{D_{\uzero}\Phi(\uzero,.).\F(.)}{\l}\leq2^{\l+\ta\nu}\m_0\CSob \NH{\F}{\l}.$$
\end{lem}
\textbf{Proof.} Let $\F$ be a function of $\Sob^\l(\R/\T\Z,\Ezero)$. Then, with the aid of lemmas \ref{propnorm1} and \ref{propnorm4}, we get
\begin{eqnarray}
\NH{D_{\uzero}\Phi(\uzero,.).\F(.)}{\l}&\leq&\sum_{\j=2}^{\p}\NH{D_{\uzero}\Phi_\j(\uzero,.).\F(.)}{\l}\nonumber\\
				&\leq&\sum_{\j=2}^{\p}\NHP{D_{\X}\Phi_\j(\X,.).\F(.)}{\j-1}\NR{\uzero}^{\j-1}\nonumber\\
				&\leq& \sum_{\j=2}^{\p}\CSob\sqrt{\m_0}\j\NHP{\Phi_\j(\X,.)}{\j}\NHP{\F(.)}{0}\NR{\uzero}^{\j-1}\nonumber\\
				&=& \sum_{\j=2}^{\p}\CSob\sqrt{\m_0}\j\pphi_\j\delta^{\j-1}\NH{\F}{\l}\nonumber
\end{eqnarray}
where, by lemma \ref{lemK}, we have
\begin{eqnarray}
\sum_{\j=2}^{\p}\j\pphi_\j\delta^{\j-1}&\leq&\sum_{\j=2}^{\p}\j\sqrt{\m_0}(\j!)^{\l+1+\ta\nu}K^{\j-1}\delta^{\j-1}\nonumber\\
			&\leq&\sqrt{\m_0}\sum_{\j=2}^{\p}\p(2\p^{\j-2})^{\l+1+\ta\nu}(K\delta)^{\j-1}\nonumber\\
			&\leq&\sqrt{\m_0}\p\left(\frac{2}{\p}\right)^{\l+1+\ta\nu}\sum_{\j=2}^{\p}(\p^{\l+1+\ta\nu}K\delta)^{\j-1}.\nonumber
\end{eqnarray}
Then, since $\p=\p_{opt}$ and $\delta<\delta_0$, we have
$$1\leq\p \quad\text{and}\quad \p^{\l+1+\ta\nu}K\delta\leq\frac{1}{2}.$$
Hence, 
$$\sum_{\j=2}^{\p}\j\pphi_\j\delta^{\j-1}\leq\sqrt{\m_0}2^{\l+\ta\nu};$$
and finally,
$$\NH{D_{\uzero}\Phi(\uzero,.).\F(.)}{\l}\leq2^{\l+\ta\nu}\m_0\CSob \NH{\F}{\l}.$$\cqfd

Then, to compute an estimate for $\V^1$, it remains to find an upper bound for
$$\NH{\V(\uzero+\vun+\Phi(\uzero,.),.)-\V(\uzero+\Phi(\uzero,.),.)}{\l}.$$
We show the following
\begin{lem}\label{lem2}
There exists $\M_1$ such that for all $\NR{\uzero},\NR{\vun}\leq\delta_0$, 
$$\NH{\V(\uzero+\vun+\Phi(\uzero,.),.)-\V(\uzero+\Phi(\uzero,.),.)}{\l}\leq\M_1\NR{\vun}\left(\NR{\uzero}+\NR{\vun}\right).$$
\end{lem}
\textbf{Proof.}
Since $\V$ belongs to $\mathcal{A}(\Nu,\Sob^\l(\R/\T\Z, \R^\m))$, and since $(A_1)\Rightarrow(A_3)$ holds in lemma \ref{propnorm3}, then for all $f,g$ in $\Sob^\l(\R/\T\Z,\R^\m)$
\begin{eqnarray}
&&\hspace{-5ex}\NH{\V_\q(.)[f+g,\cdots,f+g]-\V_\q(.)[f,\cdots,f]}{\l}\nonumber\\
&\leq&\NH{\V_\q(.)[f+g,f+g,\cdots,f+g]-\V_\q(.)[f,f+g,\cdots,f+g]}{\l}\nonumber\\
	&&+\cdots+\NH{\V_\q(.)[f,\cdots,f,f+g]-\V_\q(.)[f,\cdots,f,f]}{\l}\nonumber\\
&=&\NH{\V_\q(.)[g,f+g,\cdots,f+g]}{\l}+\NH{\V_\q(.)[f,g,f+g,\cdots,f+g]}{\l}\nonumber\\
&&+\cdots+\NH{\V_\q(.)[f,\cdots,f,g,f+g]}{\l}+\NH{\V_\q(.)[f,f,\cdots,f,g]}{\l}\nonumber\\
&\leq& \q\c\left(\frac{\CSob\sqrt{\m}}{\rh}\right)^{\q}\NH{g}{\l}\left(\NH{f}{\l}+\NH{g}{\l}\right)^{\q-1}.\nonumber	
\end{eqnarray}
Hence, with $f(t)=\uzero+\Phi(\uzero,t)=\PPhi(\uzero,t)$ and $g(t)=\vun$, we get, using lemma \ref{popt},
\begin{eqnarray}
&&\hspace{-5ex}\NH{\V(\uzero+\vun+\Phi(\uzero,.),.)-\V(\uzero+\Phi(\uzero,.),.)}{\l}\nonumber\\
&\leq&\sum_{\q=2}^{+\infty}\q\c\left(\frac{\CSob\sqrt{\m}}{\rh}\right)^{\q}\NH{\vun}{\l}\left(\NH{\uzero+\Phi(\uzero,.)}{\l}+\NH{\vun}{\l}\right)^{\q-1}\nonumber\\
&\leq&\NR{\vun}\c\left(\frac{\CSob\sqrt{\m}}{\rh}\right)^{2}\left(2\NR{\uzero}\sqrt{\m_0}+\NR{\vun}\right)\sum_{\q=2}^{+\infty}\q\left(\frac{\CSob\sqrt{\m}}{\rh}(2\delta\sqrt{\m_0}+\NR{\vun})\right)^{\q-2}.\nonumber
\end{eqnarray}
So, setting
$$\M_1:=2\c\sqrt{\m_0}\left(\frac{\CSob\sqrt{\m}}{\rh}\right)^{2}\sum_{\k=0}^{+\infty}(\k+2)\left(\frac{\CSob\sqrt{\m}}{\rh}(3\delta_0\sqrt{\m_0})\right)^{\k},$$
in which the sum converges because of the value of $\delta_0$ chosen in (\ref{defdelta0}), we get that for $\NR{\vun}\leq\delta_0\sqrt{\m_0}$,
$$\NH{\V(\uzero+\vun+\Phi(\uzero,.),.)-\V(\uzero+\Phi(\uzero,.),.)}{\l}\leq\M_1\NR{\vun}\left(\NR{\uzero}+\NR{\vun}\right)$$
holds. \cqfd

\vspace{1ex}

Thus finally, combining lemmas \ref{lem1} and \ref{lem2}, and setting
$$\M_0:=(1+||P_0||+2^{\l+\ta\nu}\m_0||P_0||)\M_1,$$
we obtain that the inequality (\ref{majvun}) holds.

\cqfd

\section{Proof of theorem \ref{thm2}}\label{proof2}

This part is entirely devoted to the proof of theorem \ref{thm2}. This proof begins by showing technical lemmas, and then we follow the same strategy of proof as that of Iooss and Lombardi (2005). Moreover, we use the same norms as those defined in part \ref{defnorm}. 

\subsection{Notations, strategy of construction for $\Phi$ and $\Nrond$}\label{constr}

First, here also, we fix $\delta$ and $\p$ and look for $\Phi$ and $\Nrond$ in 

$\mathcal{P}_{\p}(\R^\m,\Sob^{\l}(\R/\T\Z,\R^\m))$ of the form
\begin{equation}
\Phi(\X,t):=\sum_{\n=2}^{\p}\Phi_\n(\X,t),\quad \Nrond(\X,t):=\sum_{\n=2}^{\p}\Nrond_\n(\X,t),\nonumber
\end{equation}
with $\Phi_\n, \NN_\n$ in $\Hom^\n(\R^\m,\Sob^\l(\R/\T\Z,\R^\m))$ (space of homogeneous polynomials of degree $\p$). Then, one can check that the change of variables 
$$u=\y+\Phi(\y,t)$$
transforms our system (\ref{syst}) into (\ref{systNF}) if and only if
\begin{equation}\label{NFeq}
(\partial_t\hspace{-0,3ex}+\hspace{-0,3ex}\B_L)\Phi(\X,t)+(\Id\hspace{-0,3ex}+\hspace{-0,3ex}D_\X\Phi(\X,t))(\Nrond(\X,t)\hspace{-0,3ex}+\hspace{-0,3ex}\Reste(\X,t))=\V(\X\hspace{-0,5ex}+\hspace{-0,2ex}\Phi(\X,t))
\end{equation}
holds, where we have set:
$$\B_L\Phi(\X,t):=D_\X\Phi(\X,t).L\X-L\Phi(\X,t).$$
Since we look for $\Phi$ and $\Nrond$ in $\mathcal{P}_{\p}(\R^\m,\Sob^{\l}(\R/\T\Z,\R^\m))$, and $\Reste$ of order more than $\p$ in $u$, then equation (\ref{NFeq}) is equivalent to the following system
\begin{equation}\label{NFeq1}
(\partial_t+\B_L)\Phi(\X,t)\hspace{-0,5ex}+\hspace{-0,5ex}\Nrond(\X,t)\hspace{-0,5ex}=\hspace{-0,5ex}\Pi_\p(\V(\X+\Phi(\X,t))\hspace{-0,5ex}-\hspace{-0,5ex}D_\X\Phi(\X,t).\Nrond(\X,t)),
\end{equation}
\begin{equation}\label{NFeq2}
(\Id\hspace{-0,5ex}+\hspace{-0,5ex}D_\X\Phi(\X,t))\Reste(\X,t)\hspace{-0,5ex}=\hspace{-0,5ex}(\Id\hspace{-0,5ex}-\hspace{-0,5ex}\Pi_\p)(\V(\X\hspace{-0,5ex}+\hspace{-0,5ex}\Phi(\X,t))\hspace{-0,5ex}-\hspace{-0,5ex}D_\X\Phi(\X,t).\Nrond(\X,t))
\end{equation}
where $\Pi_\p$ stands for the projection on $\mathcal{P}_{\p}(\R^\m,\Sob^{\l}(\R/\T\Z,\R^\m))$. 
We begin by solving (\ref{NFeq1}). Then, (\ref{NFeq2}) will be the definition of $\Reste$ if one shows that $(\Id+D_\X\Phi(\X,t))$ is invertible (and it is; see section \ref{fin}) .
Here again, we project equation (\ref{NFeq1}) on the spaces $\Hom^\n(\R^\m,\Sob^\l(\R/\T\Z,\R^\m))$. Denoting by $\pi_\n$ this projection, we obtain
\begin{equation}\label{NFprojn}
(\partial_t\hspace{-0,3ex}+\hspace{-0,3ex}\B_L)\Phi_\n(\X,t)\hspace{-0,3ex}+\hspace{-0,3ex}\Nrond_\n(\X,t)\hspace{-0,3ex}=\hspace{-0,3ex}\pi_\n(\V(\X\hspace{-0,5ex}+\hspace{-0,3ex}\Phi(\X,t))\hspace{-0,3ex}-\hspace{-0,3ex}D_\X\Phi(\X,t).\Nrond(\X,t)).
\end{equation}
Expanding the right hand side of (\ref{NFprojn}) in power series, one can observe that, since
$$\V(0,.)=0=D_u\V(0,.),$$
(because of (\ref{nonlin})), this right hand side of (\ref{NFprojn}) only depends on $\Phi_2,\cdots,\Phi_{\n-1}$, $\Nrond_2,\cdots,\Nrond_{\n-1}$ and $\X,t$. Hence, (\ref{NFprojn}) should enables us to construct the $\Phi_\n$ and $\Nrond_\n$ by induction. We only miss the lemma proved in the following subsection to be completly convinced.

\subsection{Affine equation on $\Hom^\n(\R^\m,\Sob^\l(\R/\T\Z,\R^\m))$}\label{lem}

\begin{lem}\label{lemNF}
Let $\F_\n$ be a polynomial in $\Hom^\n(\R^\m,\Sob^\l(\R/\T\Z,\R^\m))$. There exist $\Phi_\n$ in $\Hom^\n(\R^\m,\Sob^{\l+1}(\R/\T\Z,\R^\m))$ and $\Nrond_\n$ in $\Hom^\n(\R^\m,\Sob^\l(\R/\T\Z,\R^\m))$ such that
$$(\partial_t+\B_L)\Phi_\n+\Nrond_\n=\F_\n,$$
with $\Nrond_\n$ in $\ker(-\partial_t+\B_{L^*})$, and satisfying the estimates
\begin{eqnarray}
\NHP{\Nrond_\n}{\n}&\leq&\NHP{\F_\n}{\n},\nonumber\\
\NHP{\Phi_\n}{\n}&\leq&\C_\l\n^{\l+\ta}\NHP{\F_\n}{\n}\nonumber
\end{eqnarray}
with
\begin{eqnarray}
\C_\l&:=&\max\{1,\frac{1}{\ga}\} \left(1+\frac{\T^2}{2\pi^2}(1+4\Lambda^2)\right)^{\frac{\l}{2}}\nonumber\\
\Lambda&:=&\max\limits_{1\leq\j\leq\m}|\lambda_\j|.\nonumber
\end{eqnarray}
\end{lem}

\textbf{Proof.} Here again, the key idea is that, using Fourier theory, the affine equation
\begin{equation}\label{equaffine}
(\partial_t+\B_L)\Phi_\n+\Nrond_\n=\F_\n 
\end{equation}
in $\Hom^\n(\R^\m,\Sob^\l(\R/\T\Z,\R^\m))$ is transformed into an infinity of affine equations in $\Hom^\n(\R^\m,\R^\m)$, so that we can use results of Iooss and Lombardi (2005) for each of these equations.

\begin{description}
\item[Step 1: Spliting our problem in an infinity of subproblems.]
\end{description}

Denote
\begin{eqnarray}
\Phi_\n(\X,t)&:=&\sum_{|\alpha|=\n}\phi_\alpha(t)\X^\alpha, \quad \Nrond_\n(\X,t):=\sum_{|\alpha|=\n}\Nn_\alpha(t)\X^\alpha,\nonumber\\ \F_\n(\X,t)&:=&\sum_{|\alpha|=\n}f_\alpha(t)\X^\alpha.\nonumber
\end{eqnarray}
We want all that functions to belong to $\Hom^\n(\R^\m,\Sob^\l(\R/\T\Z,\R^\m))$, with $\l\geq1$, then necessarily we can expand it in Fourier series
\begin{eqnarray}
\Phi_\n(\X,t)&=&\sum_{|\alpha|=\n}\phi_\alpha(t)\X^\alpha=\sum_{|\alpha|=\n}\bigg(\sum_{\k=-\infty}^{+\infty}\phi_\alpha^{(\k)}\E^{\I\k\frac{2\pi}{\T}t}\bigg)\X^\alpha\nonumber\\
&=&\sum_{\k=-\infty}^{+\infty}\bigg(\sum_{|\alpha|=\n}\phi_\alpha^{(\k)}\X^\alpha\bigg)\E^{\I\k\frac{2\pi}{\T}t}:=\sum_{\k=-\infty}^{+\infty}\Phi_\n^{(\k)}(\X)\E^{\I\k\frac{2\pi}{\T}t};\nonumber
\end{eqnarray}
and, with the same notations
\begin{eqnarray}
\Nrond_\n(\X,t)&:=&\sum_{\k=-\infty}^{+\infty}\Nrond_\n^{(\k)}(\X)\E^{\I\k\frac{2\pi}{\T}t},\nonumber\\
\F_\n(\X,t)&:=&\sum_{\k=-\infty}^{+\infty}\F_\n^{(\k)}(\X)\E^{\I\k\frac{2\pi}{\T}t}.\nonumber
\end{eqnarray}
Then, (\ref{equaffine}) holds if and only if
\begin{equation}\label{NFsub}
\text{for all }\k\in\Z , \quad (\I\k\frac{2\pi}{\T}+\B_L)\Phi_\n^{(\k)}(\X)+\Nrond_\n^{(\k)}=\F_\n^{(\k)}(\X)
\end{equation}
hold. 

\textbf{Now we first solve, in "Step 2", the subproblem given by (\ref{NFsub}) when we fix one $\k$, and then we show in Steps 3 to 6 that we can sum on $\k$ the Fourier series obtained.}

\begin{description}
\item[Step 2: Subproblems in $\Hom^\n(\R^\m,\R^\m)$, construction of a solution]
\end{description}

Fix any $\k$ in $\Z$. Define
\begin{eqnarray}
\Bnk:\Hom^{\n}(\R^\m,\R^\m)&\longrightarrow&\Hom^{\n}(\R^\m,\R^\m)\nonumber\\
																				\Phi(\X) &\longmapsto& (\I\k\frac{2\pi}{\T}+\B_L)\Phi(\X). \label{defAk}
\end{eqnarray}
For a given $\F_\n^{(\k)}(\X)$ in $\Hom^{\n}(\R^\m,\R^\m)$, our aim is to find $\Nrond_\n^{(\k)}(\X)$ and $\Phi_\n^{(\k)}(\X)$ such that 
\begin{equation}\label{equak}
\Bk\Phi_\n^{(\k)}(\X)+\Nrond_\n^{(\k)}=\F_\n^{(\k)}(\X).
\end{equation}
Here unfortunately, $\Bnk$ is not necessarily invertible. Hence, our strategy is to chose $\Nrond_\n^{(\k)}(\X)$ such that $\F_\n^{(\k)}(\X)-\Nrond_\n^{(\k)}(\X)$ belongs to $\Img\Bnk$, and then to take $\Phi_\n^{(\k)}(\X)$ such that $\Bk\Phi_\n^{(\k)}(\X)=\F_\n^{(\k)}(\X)-\Nrond_\n^{(\k)}(\X)$. So it is sufficient to chose a supplementary space of $\Img\Bnk$ in $\Hom^{\n}(\R^\m,\R^\m)$: indeed, then (\ref{equak}) would simply be the spliting of $\F_\n^{(\k)}(\X)$ on our decomposition of $\Hom^{\n}(\R^\m,\R^\m)$ in supplementary spaces. Making a choice of supplementary space determines in fact the normal form characterization. Here, to obtain the criteria (\ref{critereNF}) and the estimates given in the theorem, we make the following choice.

As done in Iooss and Lombardi (2005), we introduce a scalar product in $\Hom^\n(\R^\m,\R^\m)$, and then write
\begin{eqnarray}
&&\hspace{-4ex}\Hom^\n(\R^\m,\R^\m)=\Img\Bnk\hspace{-0,5ex}\stackrel{\bot}{\oplus}(\Img\Bnk)^{\bot}\hspace{-0,5ex}=\Img\Bnk\stackrel{\bot}{\oplus}\ker(\Bnk)^*,\nonumber\\
&&\hspace{-4ex}\Hom^\n(\R^\m,\R^\m)=\ker\Bnk\hspace{-0,5ex}\stackrel{\bot}{\oplus}(\ker\Bnk)^{\bot}\hspace{-0,5ex}=\ker\Bnk\stackrel{\bot}{\oplus}\Img(\Bnk)^*\hspace{-1ex}.\nonumber
\end{eqnarray}
Thus, the map $\Bnk$ reads
$$\Bnk:\ker\Bnk\stackrel{\bot}{\oplus}\Img(\Bnk)^*\longrightarrow\Img\Bnk\stackrel{\bot}{\oplus}\ker(\Bnk)^*,$$
and we can define an invertible map $\widetilde{\Bnk}$, whose inverse is called the pseudo-inverse of $\Bnk$
$$\widetilde{\Bnk}:\Img(\Bnk)^*\longrightarrow\Img\Bnk.$$
We denote by $\pi_\n^\k$ the orthogonal projection on $\ker(\Bnk)^*$, and chose
\begin{eqnarray}
\Nrond_\n^{(\k)}(\X)&:=&\pi_\n^\k(\F_\n^{(\k)}(\X)),\label{defPhi}\\
\Phi_\n^{(\k)}(\X)&:=&(\widetilde{\Bnk})^{-1}\left((\Id-\pi_\n^{\k})(\F_\n^{(\k)})\right)\label{defNrond}.
\end{eqnarray}

So, let us chose an apropriate scalar product. We chose the same scalar product as that introduced in 2.1 of Iooss and Lombardi (2005). Namely, for any pair of polynomials $P$ and $P'$, define
\begin{equation}\label{defscal}
\scal{P}{P'}=\overline{P}(\partial_\X)P'(\X)|_{\X=0}.
\end{equation}
Moreover, the norm associated with this scalar product is $\NP{.}{\n}$, introduced in section \ref{defnorm}.

\begin{description}
\item[Step 3: Subproblems in $\Hom^\n(\R^\m,\R^\m)$, upper bounds]
\end{description}

To sum the $\Nrond_\n^{(\k)}(\X)$ and $\Phi_\n^{(\k)}(\X)$ defined in (\ref{defPhi}), (\ref{defNrond}) , we need to compute upper bounds.
Since $\pi_\n^\k$ is an orthogonal projection, we have:
\begin{eqnarray}
\NP{\Nrond_\n^{(\k)}(\X)}{\n}&\leq&\NP{\F_\n^{(\k)}(\X))}{\n},\nonumber\\
\NP{\Phi_\n^{(\k)}(\X)}{\n}&\leq&\NNorme{(\widetilde{\Bnk})^{-1}}{2,\n}\NP{\F_\n^{(\k)}}{\n}.\nonumber
\end{eqnarray}
Moreover, lemma 2.5 of Iooss and Lombardi (2005) ensures that, since by hypothesis \ref{hyp3} $L$ is diagonizable, then $\mathcal{B}_L|_{\Hom^\n}$ is so (our $\mathcal{B}_L$ is their $\mathcal{A}_L$), with eigenvalues:
\begin{equation}\label{vpAL}
\{\scal{\a}{\lambda}-\lambda_\j,  1\leq\j\leq\m, \a\in\N^\m, |\a|=\n\}.
\end{equation}
Hence, the operator $\Bnk$ defined in (\ref{defAk}) is diagonizable, with eigenvalues
$$\{\scal{\a}{\lambda}+\I\k\frac{2\pi}{\T}-\lambda_\j,  1\leq\j\leq\m, \a\in\N^\m, |\a|=\n\}.$$
Then, using (\ref{diophantian2}) of hypothesis \ref{hyp3}, we get
\begin{equation}
\NNorme{(\widetilde{\Bnk})^{-1}}{2,\n}\leq\max\limits_{\dindice{1\leq\j\leq\m}{\a\in\N^\m, |\a|=\n}}\frac{1}{|<\a,\lambda>+\I\k\frac{2\pi}{\T}-\lambda_\j|}\leq\frac{(\n+|\k|)^\ta}{\ga}.\nonumber
\end{equation}
So we finally have the upper bounds:
\begin{eqnarray}
\NP{\Nrond_\n^{(\k)}(\X)}{\n}&\leq&\NP{\F_\n^{(\k)}(\X)}{\n},\label{NFmajNrond}\\
\NP{\Phi_\n^{(\k)}(\X)}{\n}&\leq&\frac{(\n+|\k|)^\ta}{\ga}\NP{\F_\n^{(\k)}(\X)}{\n}.\label{NFmajPhi}
\end{eqnarray}

\begin{description}
\item[Step 4: $\Nrond_\n$ and $\Phi_\n$ are well-defined, $\Nrond_\n\in\Hom^\n(\R^\m,\Sob^\l(\R/\T\Z,\R^\m))$]
\end{description}
We now set
\begin{eqnarray}
\hspace{-6ex}\Nrond_\n(\X,t)\hspace{-2ex}&:=&\hspace{-2ex}\sum_{\k=-\infty}^{+\infty}\Nrond_\n^{(\k)}(\X)\E^{\I\k\frac{2\pi}{\T}t}=\hspace{-1ex}\sum_{\k=-\infty}^{+\infty}\hspace{-1ex}\pi_\n^\k(\F_\n^{(\k)}(\X))\E^{\I\k\frac{2\pi}{\T}t},\label{defNrondn}\\
\hspace{-6ex}\Phi_\n(\X,t)\hspace{-2ex}&:=&\hspace{-3ex}\sum_{\k=-\infty}^{+\infty}\Phi_\n^{(\k)}(\X)\E^{\I\k\frac{2\pi}{\T}t}=\hspace{-2ex}\sum_{\k=-\infty}^{+\infty}\hspace{-1ex}\widetilde{\Bk}^{-1}\hspace{-1ex}\left((\Id-\pi_\n^{\k})(\F_\n^{(\k)}(\X))\right)\E^{\I\k\frac{2\pi}{\T}t}\label{defPhin}
\end{eqnarray}
and we want to prove that they satisfy all the properties required in lemma \ref{lemNF}. We begin by showing that they are well-defined, verifying that $\Nrond_\n$ is in $\Hom^\n(\R^\m,\Sob^\l(\R/\T\Z,\R^\m))$ and $\Phi_\n$ in $\Hom^\n(\R^\m,\Sob^0(\R/\T\Z,\R^\m))$, with the norms defined in section \ref{defnorm}. 

Using the upper bounds (\ref{NFmajNrond}) and (\ref{NFmajPhi}) computed in step 3, and the inequality $\ta\leq\l$ assumed by hypothesis \ref{hyp3}, we get:
\begin{eqnarray}
\NHP{\Nrond_\n(\X,t)}{\n}^2&=&\sum_{\k=-\infty}^{+\infty}(1+\k^2)^\l\NP{\Nrond_\n^{(\k)}(\X)}{\n}^2\nonumber\\
					&\leq&\sum_{\k=-\infty}^{+\infty}(1+\k^2)^\l\NP{\F_\n^{(\k)}(\X)}{\n}^2=\NHP{\F_\n(\X,t)}{\n}^2<+\infty;\nonumber
\end{eqnarray}
\begin{eqnarray}
\NHP{\Phi_\n(\X,t)}{\n,\Sob^0}^2&=&\sum_{\k=-\infty}^{+\infty}\NP{\Phi_\n^{(\k)}(\X)}{\n}^2\leq\sum_{\k=-\infty}^{+\infty}\frac{(\n+|\k|)^{2\ta}}{\ga^2}\NP{\F_\n^{(\k)}(\X)}{\n}^2\nonumber\\
					&\leq&\frac{\n^{2\ta}}{\ga^2}\sum_{\k=-\infty}^{+\infty}(1+\k^2)^\ta\NP{\F_\n^{(\k)}(\X)}{\n}^2\nonumber\\
					&\leq&\frac{\n^{2\ta}}{\ga^2}\sum_{\k=-\infty}^{+\infty}(1+\k^2)^\l\NP{\F_\n^{(\k)}(\X)}{\n}^2\nonumber\\
\NHP{\Phi_\n(\X,t)}{\n,\Sob^0}^2&\leq&\frac{\n^{2\ta}}{\ga^2}\NHP{\F_\n(\X,t)}{\n}^2<+\infty.\label{equainit}	
\end{eqnarray}

\begin{description}
\item[Step 5: $\Phi_\n$ belongs to $\Hom^\n(\R^\m,\Sob^{\l+1}(\R/\T\Z,\R^\m))$]
\end{description}

Let us prove by induction that if $0\leq\j\leq\l+1$ then
\begin{equation}\label{NFrec}
\NHP{\Phi_\n(\X,t)}{\n,\Sob^\j}\leq\C_\j\n^{\j+\ta}\NHP{\F_\n(\X,t)}{\n},
\end{equation}
where
$$\C_\j:=\max\{1,\frac{1}{\ga}\} \left(1+\frac{\T^2}{2\pi^2}(1+4\Lambda^2)\right)^{\frac{\j}{2}}.$$

Inequality (\ref{equainit}) ensures that (\ref{NFrec}) holds for $\j=0$. Asume that (\ref{NFrec}) holds for one $\j$, $0\leq\j\leq\l$.

For all $\k$ in $\Z$, since $\Phi_\n$ is defined by (\ref{defPhin}), we have
$$\Bk\Phi_\n^{(\k)}=(\Id-\pi_\n^\k)\F_\n^{(\k)};$$
and hence
$$\I\k\frac{2\pi}{\T}\Phi_\n^{(\k)}=(\Id-\pi_\n^\k)\F_\n^{(\k)}-\B_L\Phi_\n^{(\k)}.$$
Thus
$$(1+\k^2)\NP{\Phi_\n^{(\k)}}{\n}^2=\left(\frac{\T}{2\pi}\right)^2\NP{(\Id-\pi_\n^\k)\F_\n^{(\k)}-\B_L\Phi_\n^{(\k)}}{\n}^2+\NP{\Phi_\n^{(\k)}}{\n}^2;$$
where the eigenvalues of $\B_L|_{\Hom^\n}$ are given in (\ref{vpAL}). Then we get:
\begin{eqnarray}
\NP{\B_L\Phi_\n^{(\k)}}{\n}&\leq&\max\limits_{\dindice{1\leq\j\leq\m}{\a\in\N^\m, |\a|=\n}}\{|\left\langle \a,\lambda\right\rangle-\lambda_\j|\}\NP{\Phi_\n^{(\k)}}{\n}\nonumber\\
&\leq&(\n+1)\max\limits_{1\leq\j\leq\m}|\lambda_\j|\NP{\Phi_\n^{(\k)}}{\n}\hspace{-1ex}=(\n+1)\Lambda\NP{\Phi_\n^{(\k)}}{\n}.\nonumber
\end{eqnarray}
So
\begin{eqnarray}
(1+\k^2)\NP{\Phi_\n^{(\k)}}{\n}^2&\leq&\left(\frac{\T}{2\pi}\right)^2\left(\NP{\F_\n^{(\k)}}{\n}+(\n+1)\Lambda\NP{\Phi_\n^{(\k)}}{\n}\right)^2+\NP{\Phi_\n^{(\k)}}{\n}^2\nonumber\\
	&\leq&2\left(\frac{\T}{2\pi}\right)^2\hspace{-1ex}\left(\NP{\F_\n^{(\k)}}{\n}^2+(\n+1)^2\Lambda^2\NP{\Phi_\n^{(\k)}}{\n}^2\right)\hspace{-0,5ex}+\NP{\Phi_\n^{(\k)}}{\n}^2\nonumber\\
	&\leq&\frac{\T^2}{2\pi^2}\NP{\F_\n^{(\k)}}{\n}^2+\left(1+\frac{\Lambda^2\T^2(\n+1)^2}{2\pi^2}\right)\NP{\Phi_\n^{(\k)}}{\n}^2;\nonumber
\end{eqnarray}
and then
\begin{eqnarray}
&&\hspace{-7ex}\NHP{\Phi_\n}{\n,\Sob^{\j+1}}^2=\sum_{\k=-\infty}^{+\infty}(1+\k^2)^{\j+1}\NP{\Phi_\n^{(\k)}}{\n}^2\nonumber\\
&\leq&\sum_{\k=-\infty}^{+\infty}(1+\k^2)^{\j}\left(\frac{\T^2}{2\pi^2}\NP{\F_\n^{(\k)}}{\n}^2+\big(1+\frac{\Lambda^2\T^2(\n+1)^2}{2\pi^2}\big)\NP{\Phi_\n^{(\k)}}{\n}^2\right)\nonumber\\
&\leq&\frac{\T^2}{2\pi^2}\NHP{\F_\n}{\n}^2+\big(1+\frac{\Lambda^2\T^2(\n+1)^2}{2\pi^2}\big)\NHP{\Phi_\n}{\n,\Sob^\j}^2\nonumber\\
&\leq&\left(\frac{\T^2}{2\pi^2}+\big(1+\frac{\Lambda^2\T^2(\n+1)^2}{2\pi^2}\big)\C_\j^2\n^{2(\j+\ta)}\right)\NHP{\F_\n}{\n}^2\nonumber\\
&\leq& \C_{\j+1}^2\n^{2(\j+1+\ta)}\NHP{\F_\n}{\n}^2.\nonumber
\end{eqnarray}
Finally, we have proved that
$$\NHP{\Phi_\n}{\n,\Sob^{\l+1}}\leq \C_{\l}\hspace{0,5ex}\n^{\l+\ta}\NHP{\F_\n}{\n}$$
holds. This shows in particular that $\Phi_\n$ is in $\Hom^\n(\R^\m,\Sob^{\l+1}(\R/\T\Z,\R^\m))$.
\begin{description}
\item[Step 6: $\Nrond_\n$ is in $\ker(-\partial_t+\B_{L^*})$]
\end{description}

Since $\pi_\n^\k$ is the projection on $\ker(\Bnk)^*$, then for all $\k$ in $\Z$,
\begin{equation}\label{Nk}
\Nrond_\n^{(\k)}(\X)=\pi_\n^\k\left(\F_\n^{(\k)}(\X)\right)\in\ker(\Bnk)^*;
\end{equation}
where
$$(\Bk)^*=-\I\k\frac{2\pi}{\T}+(\B_L)^*.$$
And the choice of scalar product (\ref{defscal}) guarantees that 
$$(\B_L)^*=\B_{L^*}$$
holds (see proof by Iooss and Adelmeyer (1992)). Then, for all $\k$ in $\Z$, (\ref{Nk}) ensures that
$$\B_{L^*}\Nrond_\n^{(\k)}(\X)=\I\k\frac{2\pi}{\T}\Nrond_\n^{(\k)}(\X).$$
Finally, since $\Nrond_\n$ is in $\Hom^\n(\R^\m,\Sob^{\l}(\R/\T\Z,\R^\m))$ with $\l\geq1$, we can sum on $\k$, and we obtain
$$\B_{L^*}\Nrond_\n(\X,t)=\partial_t\Nrond_\n(\X,t).$$\cqfd

\subsection{End of the proof (sketch)}\label{fin}
Now, combining parts \ref{constr} and \ref{lem}, we get that we have constructed $\Phi$ in $\mathcal{P}_\p(\R^\m,\Sob^{\l+1}(\R/\T\Z,\R^\m))$ and $\Nrond$ in $\mathcal{P}_\p(\R^\m,\Sob^{\l}(\R/\T\Z,\R^\m))$, with $\Nrond$ also in $\ker(-\partial_t+\B_{L^*})$, such that if one writes
\begin{equation}
\Phi(\X,t):=\sum_{\n=2}^{\p}\Phi_\n(\X,t),\quad \Nrond(\X,t):=\sum_{\n=2}^{\p}\Nrond_\n(\X,t),\nonumber
\end{equation}
then the $\Phi_\n$ and $\Nrond_\n$ satisfy
$$(\partial_t+\B_L)\Phi_\n(\X,t)+\Nrond_\n(\X,t)=\pi_\n(\V(\X+\Phi(\X,t))-D_\X\Phi(\X,t).\Nrond(\X,t));$$
and the estimates 
\begin{eqnarray}
\NHP{\Nrond_\n(\X,t)}{\n}&\leq&\NHP{\pi_\n(\V(\X+\Phi(\X,t))-D_\X\Phi(\X,t).\Nrond(\X,t))}{\n},\nonumber\\
\NHP{\Phi_\n(\X,t)}{\n}&\leq&\C_\l\n^{\l+\ta}\NHP{\pi_\n(\V(\X+\Phi(\X,t))-D_\X\Phi(\X,t).\Nrond(\X,t))}{\n}.\nonumber
\end{eqnarray}

\vspace{1ex}

First, one can check that the fact that $\Nrond$ belongs to $\ker(-\partial_t+\B_{L^*})$ guarantees that the normal form criteria (\ref{critereNF}) holds.

Then it remains to show that the remainder $\Reste$ is well-defined by (\ref{NFeq2}) and that the estimate (\ref{NFmaj}) holds. We follow the same strategy as in part 2.3 of Iooss and Lombardi (2005): our lemmas \ref{propnorm3} and \ref{propnorm4} replace their lemma 2.11, we then compute estimates for the $\nu_\n:=\NHP{\Nrond_\n}{\n}$ and $\pphi_\n:=\NHP{\Phi_\n}{\n}$, and we finally obtain a similar proof of well-definition and upper bound for $\Reste$ with slightly different constants. Namely, one can check that this way we get the following lemma instead of their last lemma 2.21.
\begin{lem}
For all $\delta>0$, for $\p=\p_{opt}(\delta)$, with $\Phi$ and $\Nrond$ constructed above, the remainder $\Reste$ is well-defined by (\ref{NFeq2}) and satisfies
$$\supl_{\NR{\y}\leq\delta}\NH{\Reste(\y,.)}{\l}\leq\M'\delta^2\E^{-\frac{\om}{\delta^\b}}.$$
where
$$\b=\frac{1}{1+\ta'}=\frac{1}{1+\l+\ta}, \quad \p_{opt}(\delta)=\left[\frac{1}{\E(\textbf{C}\delta)^\b}\right], \quad \om=\frac{1}{\E \textbf{C}^\b},$$
$$\textbf{C}=\frac{(\CSob\sqrt{\m})^3}{\rh^2}\left[\left(\frac{5}{2}\CSob^2\m+2\right)\C_\l\c+3\frac{\rh}{\CSob\sqrt{\m}}\right],$$
$$\M'=\frac{10}{9}\c \textbf{C}^2\left(\left(\mathcal{M}\sqrt{\frac{27}{8\E}}\right)^{1+\l+\ta}+(2\E)^{2(1+\l+\ta)}\right), \quad \mathcal{M}=\supl_{\p\in\N}\frac{\E^2\p!}{\p^{\p+\frac{1}{2}}\E^{-\p}}.$$
\end{lem}


\begin{thebibliography}{99}
\bibitem{Kam2}
				Arnold, V. (1978). Mathematical Methods of Classical Mechanics.	Springer-Verlag, New-York.
\bibitem{CouetteTaylor}
				Chossat, P., and Iooss, G. (1994). The Couette-Taylor problem.	Appl. Math. Sci. 102.
\bibitem{NFcritere}
				Elphick, C., Tirapegui, E., Brachet, M.E., Coullet, P., and G.Iooss. (1987). A simple global characterization for normal forms of singular vector fields.	Physica D 29, 95-127.
\bibitem{Holmes}
				Guckenheimer, J., and Holmes, P. (1983). Nonlinear oscillations, dynamical systems and bifurcations of vector fields. Applied Mathematical        Science, vol.42, Springer, Berlin, Heidelberg, New York.
\bibitem{GerardMariana}
				Haragus, M., and Iooss, G. (2010). Local bifurcations, center manifolds, and normal forms in infinite dimensional dynamical systems. Springer Verlag, to appear.
\bibitem{Gerard}
				Iooss, G., and Adelmeyer, M. (1992). Topics in Bifurcation Theory and Applications. Advanced Series in Nonlinear Dynamics, vol. 3, World        		    Scientific, Singapore.
\bibitem{Split}
				Iooss, G., and Lombardi, E. (2010). Approximate invariant manifolds up to exponentially small terms.
				J. Differential Equations 248, 14101431.
\bibitem{NFexp} 
				Iooss, G., and Lombardi, E. (2005). Polynomial normal forms with exponentially small remainder for analytic vector fields. J. Differential 			        Equations 212, 1-61.
\bibitem{IoossPeroueme}
				Iooss, G., and P\'erou\'eme, M.C. (1993). Perturbed homoclinic solutions in 1:1 resonance vector fields.	J. Differential Equations 102(1).	
\bibitem{center}
				Kelley, A. (1967). The stable, center-stable, center, center-unstable, unstable manifolds. J. Differential Equations 3, 546-570.	
\bibitem{Kam1}
				Kolmogorov, A.N. (1954). On conservation of conditionally periodic motions under small perturbation of the hamiltonian.	Dokl. Akad. Nauk, SSSR, 98:527-530.
\bibitem{Eric}
				Lombardi, E. (2000). Oscillatory Integrals and Phenomena Beyond all Algebraic Orders.	Lecture Notes in Mathematics vol.1741. Springer.
\bibitem{Neko1}
				Nekoroshev, N.N. (1977). An exponential estimate of the time of stability of nearly integrable Hamiltonian systems, I. Usp. Mat. Nauk 32 (1997) 5-66 Russ. Math. Surv. 32 (1977) 1-65.
\bibitem{Neko2}
				Nekoroshev, N.N. (1979). An exponential estimate of the time of stability of nearly integrable Hamiltonian systems, II.	Tr. Semin. Petrovsk. 5 (1979) 5-50 in: O.A. Oleineik (Ed.), Topics in Modern Mathematics, Petrovskii Semin., no.5. Consultant Bureau, New York, 1985.				
\bibitem{Poschel}
				Poschel, J. (1993). Nekoroshev Estimates for Quasi-convex Hamiltonian Systems. Math. Z. 213  187-216.
\bibitem{Physic}
				Touz\'e, C., and Amabili, M. (2006). Nonlinear normal modes for damped geometrically nonlinear systems: Application to reduced-order modelling of          harmonically forced structures. Journal of Sound an Vibration 298, 958-981.	
\end{thebibliography}
\end{document}